\documentclass[12pt]{article}
\usepackage{hyperref}
\usepackage{amssymb, amsmath, amsthm}
\usepackage[all]{xy}
\usepackage{ifthen}

\newcommand{\mr}[1]{\mathrm{#1}}
\newcommand{\mf}[1]{\mathfrak{#1}}
\newcommand{\mc}[1]{\mathcal{#1}}

\newcommand{\Z}{{\bf Z}}
\newcommand{\Q}{{\bf Q}}
\newcommand{\zp}{{\bf Z}_p}
\newcommand{\qp}{{\bf Q}_p}
\newcommand{\qbar}{\overline{\Q}}
\newcommand{\qpbar}{\overline{\qp}}
\newcommand{\et}{\text{\'et}}
\newcommand{\ord}{\mr{ord}}
\newcommand{\cts}{\mr{cts}}

\newenvironment{ilist}{\begin{list}{}{\setlength{\leftmargin}{0ex}\setlength{\topsep}{0ex}\setlength{\itemsep}{0ex}\setlength{\parskip}{0pt}\setlength{\parsep}{0pt}}}{\end{list}}
\newenvironment{ilist2}{\begin{list}{}{\setlength{\leftmargin}{0ex}\setlength{\topsep}{0ex}\setlength{\itemsep}{0ex}\setlength{\parskip}{0pt}\setlength{\parsep}{0pt}\setlength{\itemindent}{20pt}
}}{\end{list}}

\newcommand{\pr}[1]{\pageref{#1}}
\newcommand{\sr}[2]{\ifthenelse{\equal{\pr{#1}}{\pr{#2}}}{\pr{#1}}{\pr{#1}--\pr{#2}}}

\newcommand{\cotimes}{\,\hat{\otimes}_{\zp}\,}
\newcommand{\fpath}[4]{\left\{\frac{#1}{#2},\frac{#3}{#4}\right\}}

\newcommand{\I}{\mc{I}}
\newcommand{\h}{\mf{h}^*_{\mf{m}}}
\newcommand{\ltheta}{\langle \theta \rangle}

\DeclareMathOperator{\Hom}{Hom} 
\DeclareMathOperator{\End}{End} \DeclareMathOperator{\Gal}{Gal}

\newtheorem{theorem}{Theorem}[section]
\newtheorem{proposition}[theorem]{Proposition}
\newtheorem{lemma}[theorem]{Lemma}
\newtheorem{corollary}[theorem]{Corollary}

\newtheorem{conjecture}[theorem]{Conjecture}

\theoremstyle{definition}

\theoremstyle{remark}
\newtheorem*{remark}{Remark}

\newtheorem*{ack}{Acknowledgments}

\numberwithin{equation}{section}

\renewcommand{\baselinestretch}{1.2}

\begin{document}
\title{A reciprocity map and the\\ two variable $p$-adic $L$-function}
\author{Romyar T. Sharifi}
\date{}
\maketitle

\begin{abstract}
For primes $p \ge 5$, we propose a conjecture that relates the values of cup products in the Galois cohomology of the maximal unramified outside $p$ extension of a cyclotomic field on cyclotomic $p$-units to the values of $p$-adic $L$-functions of cuspidal eigenforms that satisfy mod $p$ congruences with Eisenstein series.  Passing up the cyclotomic and Hida towers, we construct an isomorphism of certain spaces that allows us to compare the value of a reciprocity map on a particular norm compatible system of $p$-units to what is essentially the two-variable $p$-adic $L$-function of Mazur and Kitagawa.
\end{abstract}

\section{Introduction}

\subsection{Background}

The principal theme of this article is that special elements in the Galois cohomology of a cyclotomic field should correspond to special elements in the
quotient of the homology group of a modular curve by an Eisenstein ideal.  The elements on the Galois side of the picture arise as cup products of units in our cyclotomic field, while the elements on the modular side arise in alternate forms of our conjecture from Manin symbols and $p$-adic $L$-values of cusp forms that satisfy congruences with Eisenstein series at primes over $p$.  
We can also understand this as a comparison between objects that interpolate these elements: the value of a reciprocity map on a particular norm compatible sequence of $p$-units and an object giving rise to a two-variable $p$-adic $L$-function, taken modulo an Eisenstein ideal.

We make these correspondences explicit via a map
from the Galois group of the maximal unramified abelian pro-$p$ extension of the cyclotomic field of all $p$-power roots of unity to the quotient by an Eisenstein ideal of the inverse limit of first \'etale cohomology groups of modular curves of $p$-power level.   
Recall that the main conjecture of Iwasawa theory tells us
 that the $p$-adic zeta function 
provides a characteristic power series for the minus part of 
latter Galois group as an Iwasawa module.
In fact, the map we construct is a modification of that found in the work of M.\ Ohta on the main conjecture \cite{ohta-ord}-\cite{ohta-comp2}, which incorporated ideas of Harder-Pink and Kurihara and the Hida-theoretic aspects of the work of Wiles into a refinement of the original proof of Mazur-Wiles.  
In one of its various guises, our conjecture asserts that this map carries inverse limits of cup product values on special cyclotomic units to universal $p$-adic $L$-values modulo an Eisenstein ideal, up to a canonical unit.

The reasons to expect such a conjecture, though numerous, are far from obvious.  The core of this article being focused on the statements of the various forms of this conjecture and the proofs of their equivalence, we take some space in this first subsection to mention a few of the theoretical reasons that we expect the conjecture to hold.  We omit technical details, deferring them for the most part to future work.

Initial evidence for our conjecture can be seen in relation to the main conjecture for modular forms.   
In fact, we can show that cup products control the Selmer groups of certain reducible representations, such as the residual representations attached to newforms that satisfy mod $p$ congruences with Eisenstein series.  More precisely, under weak assumptions, such a Selmer group will be given as the quotient of an eigenspace of a cyclotomic class group modulo $p$ by the subgroup generated by a cup product of cyclotomic $p$-units.   On the other hand, $p$-adic $L$-values of such newforms are expected to control the structure of these Selmer groups by the main conjecture of Iwasawa theory for modular forms \cite[p.\ 291]{greenberg}.  That is, the main conjecture leads us to expect agreement between these cup products 
and the mod $p$ reductions of the $p$-adic $L$-values of these newforms inside the proper choice of lattice.

One can think of our conjecture as related to the main conjecture for modular forms, modulo an Eisenstein ideal, 
in a quite similar manner to that in which the classical main conjecture relates to Iwasawa's construction of the $p$-adic zeta function out of cyclotomic $p$-units.  Iwasawa's theorem provides an explicit map from the group of norm compatible sequences in the $p$-completions of the multiplicative groups of the $p$-adic fields of $p$-power roots of unity to the Iwasawa algebra that sends a compatible sequence of one minus $p$-power roots of unity to the $p$-adic zeta function  \cite{iwasawa}.  In our conjecture, the two variable $p$-adic $L$-function modulo an Eisenstein ideal is constructed out of a reciprocity map applied to the same sequence of cyclotomic $p$-units, or more loosely, out of cup products of cyclotomic $p$-units.

We remark that Fukaya proved a direct analogue of Iwasawa's theorem in the modular setting, constructing a certain two-variable $p$-adic $L$-function out of the Beilinson elements that appear in Kato's Euler system \cite{fukaya}.  In fact, Kato constructed maps that yield a comparison between these Beilinson elements, which are cup products of Siegel units, and $L$-values of cusp forms \cite{kato}.  The connection with our elements is seen in the fact that Siegel units specialize to cyclotomic $p$-units at cusps.  Fukaya constructed her modular two-variable $p$-adic $L$-function via a map arising from Coleman power series.  Although this map is defined entirely differently to ours, this nonetheless strongly suggests the existence of a direct correspondence of the sort we conjecture.

We feel obliged to emphasize, at this point, that the map that we use arises in a specific manner from the action of Galois on modular curves, which makes the conjecture considerably more delicate than a simple correspondence.  It is natural to ask why such a map should be expected to provide our comparison.  At present, the most convincing evidence we have of this is a proof of a particular specialization of the conjecture.  That is, one can derive from \cite[Theorem 5.2]{me-paireis} that our map takes a particular value of the cup product to a universal $p$-adic $L$-value at the trivial character under the assumption that $p$ does not divide a certain Bernoulli number, up to a given canonical unit.  We describe this just as briefly but more concretely in the next subsection.  It was this result that convinced us to look at the map we construct here.
That the values on cup products of this consequential map should have prior arithmetic interest in and of themselves is perhaps the most remarkable aspect of our conjectures.

\subsection{A special case}

We first describe a special but fundamental case.  Set $F = \Q(\mu_p)$ for an irregular prime $p$, and consider the $p$-completion $\mc{E}_F$ of the $p$-units in $F$.  The cup product in the Galois cohomology of the maximal unramified outside $p$ extension of $F$ defines a pairing
$$
	(\,\cdot\,,\,\cdot\,) \colon 
	\mc{E}_F \times \mc{E}_F \to A_F \otimes \mu_p,
$$
where $A_F$ denotes the $p$-part of the class group of $F$.  This pairing was studied in detail in \cite{mcs}.  We fix a complex embedding $\iota$ of $\qbar$ and thereby a  $p$th root of unity $\zeta_p = \iota^{-1}(e^{2\pi\sqrt{-1}/p})$.   Let $\omega$ denote the $p$-adic Teichm\"uller character.  For odd $t \in \Z$, define 
$$
	\alpha_t = \prod_{i=1}^{p-1} (1-\zeta_p^i)^{\omega(i)^{t-1}} \in \mc{E}_F.
$$
We consider the values $(\alpha_t,\alpha_{k-t})$ for odd integers $t$ and even integers $k$.
Such a value can be nontrivial only if the $\omega^{1-k}$-eigenspace of $A_F$ is nontrivial, which is to say, only if $p$ divides the generalized Bernoulli number $B_{1,\omega^{k-1}}$.  We fix such a $k$. 

Suppose we are given a newform $f$ of weight $2$, level $p$,  and character $\omega^{k-2}$ (coefficients in $\qpbar$) that satisfies a congruence with the normalized Eisenstein series with $l$th eigenvalue $1+\omega^{k-2}(l)l$ for odd primes $l \neq p$.  
Inside the $p$-adic representation attached to $f$ is a choice of lattice that corresponds to the first \'etale cohomology group of the closed modular curve $X_1(p)$ over $\qbar$.  The action of Galois on the resulting residual representation $T_f$ gives rise directly to a map
$$
	A_F^- \to \Hom_{\zp}(T_f^+,T_f^-),
$$
where $T_f^{\pm}$ are the $(\pm 1)$-eigenspaces of $T_f$ under complex conjugation.  We find a generator of $T_f^+$, canonical up to $\iota$, and therefore we obtain a map
$$
	\phi_f \colon A_F^- \otimes \mu_p \to T_f^- \otimes \mu_p.
$$ 

We conjecture that the map $\phi_f$ takes values of the cup product on our special cyclotomic units to the images of $p$-adic $L$-values of $f$ in $T_f^- \otimes \mu_p$.
In particular, the space $T_f^-(1)$ may be thought of as a space in which the $p$-adic $L$-values $L_p(f,\chi,s)$ naturally lie, for $\chi$ an even character and $s \in \zp$.  Let us denote the image of $L_p(f,\chi,s)$ in $T_f^- \otimes \mu_p$ by $\overline{L_p(f,\chi,s)}$.  In this setting, our conjectures state that
\begin{equation} \label{simpleconj}
	\phi_f((\alpha_t,\alpha_{k-t})) = c_{p,k}\cdot \overline{L_p(f,\omega^{t-1},1)}
\end{equation}
for some $c_{p,k} \in (\Z/p\Z)^{\times}$ independent of $t$ and $f$.

The first theoretical piece of evidence for this conjecture may be derived from \cite[Theorem 5.2]{me-paireis}.  It implies that \eqref{simpleconj} holds for $t = 1$ for some $c_{p,k} \in (\Z/p\Z)^{\times}$, under the assumption that $p$ does not also divide $B_{1,\omega^{1-k}}$.  Moreover, one can show that the value $(\alpha_t,\alpha_{k-t})$ is zero only if the Selmer group over $\Q$ of the Tate twist $T_f(t)$ of $T_f$ is nonzero under certain mild assumptions.  On the other hand, that $\overline{L_p(f,\omega^{t-1},1)}$ is zero only if the same  
Selmer group is nonzero would follow in this case from the main conjecture of Iwasawa theory for modular forms.  We intend to explore an Iwasawa-theoretic generalization of this in forthcoming work.

\subsection{Summary of the conjectures}

Let us now turn to the general setting and give a condensed but nearly precise overview of the objects to be studied in our conjectures.  Choose a prime 
$p \ge 5$ and a positive integer $N$ prime to $p$ with $p$ not dividing 
the number $\varphi(N)$ of positive integers relatively prime and less than or equal
to $N$.
The different versions of the conjecture can roughly be stated as giving, respectively, the following correspondences between to-be-defined objects:
\begin{eqnarray} 	
	(1-\zeta_{Np^r}^i,1-\zeta_{Np^r}^j)_{F_r,S}^{\circ} &\longleftrightarrow&
	\bar{\xi}_r(i:j) \label{corres1} \\
	\Psi_K^{\circ}(1-\zeta) &\longleftrightarrow& \overline{\mc{L}_N^{\star}}
	\label{corres2}\\
	\alpha_t^{\psi} \cup \alpha_{k-t}^{\theta\psi^{-1}\omega^{-1}}
	&\longleftrightarrow& \overline{L_p(\xi,\omega\theta,k,\psi,t)}.
	\label{corres3}
\end{eqnarray}
In the rest of this introduction, we first sketch the definition of the objects on the Galois (left) side of the picture, followed by the objects on the modular (right) side, and finish by describing the maps yielding the correspondences.

Let $K = \Q(\mu_{Np^{\infty}})$.  A fixed choice of complex embedding affords us norm compatible choices $\zeta_{Np^r}$ of primitive $Np^r$th roots of unity in the fields $F_r = \Q(\mu_{Np^r})$ for $r \ge 1$.  We let $S$ denote the set of primes over $Np$ and any real places of any given number field, and we let $G_{F_r,S}$ denote the Galois group of the maximal unramified outside $S$ extension of $F_r$.  We then form the cup product
$$
	H^1_{\cts}(G_{F_r,S},\zp(1))^{\otimes 2} \to 
	H^2_{\cts}(G_{F_r,S},\zp(2)).
$$
We use $(\,\cdot\,,\,\cdot\,)_{F_r,S}^{\circ}$ to denote the projection of the resulting pairing on $S$-units of $F_r$ to the sum of odd, primitive eigenspaces of the second cohomology group under a twist by $\zp(-1)$ of the standard action of $\Gal(F_1/\Q) \cong (\Z/Np\Z)^{\times}$ (see Section \ref{cpms}).  In \eqref{corres1}, we then consider values $(1-\zeta_{Np^r}^i,1-\zeta_{Np^r}^j)_{F_r,S}^{\circ}$ of this pairing for $i, j \in \Z$ nonzero modulo $Np^r$ with $(i,j,Np)=1$.

Inverse limits of these cup product pairings up the cyclotomic tower allow us to define a certain reciprocity map $\Psi_K^{\circ}$ on norm compatible sequences of $p$-units in intermediate extensions of $K/F$.  Put another way, we consider an exact sequence 
$$
	1 \to \mf{X}_K \to \mc{T} \to \zp \to 0
$$
of $\zp[[G_{K,S}]]$-modules, where
$\mf{X}_K$ is the maximal abelian pro-$p$ quotient of $G_{K,S}$, on which $G_{K,S}$
acts trivially, and
$\mc{T}$ is determined by the cocycle that is the projection map from $G_{K,S}$ to
$\mf{X}_K$.  It yields a long exact sequence among inverse limits under corestriction of cohomology groups of the $G_{F_r,S}$, in particular a coboundary map (see Section \ref{recmap})
$$
	\Psi_K \colon 
	\lim_{\leftarrow} H^1_{\cts}(G_{F_r,S},\zp(1)) \to 
	\lim_{\leftarrow} H^2_{\cts}(G_{F_r,S},\zp(1)) \otimes_{\zp} \mf{X}_K
$$  
after twisting by $\zp(1)$.  The odd, primitive part of the latter inverse limit is isomorphic to the odd, primitive part $X_K^{\circ}$ of the maximal unramified quotient $X_K$ of $\mf{X}_K$.  Then $\Psi_K^{\circ}$ is given by composing with projection to $X_K^{\circ} \otimes_{\zp} \mf{X}_K^-$, where $\mf{X}_K^-$ denotes the odd part of $\mf{X}_K$.  We are interested in \eqref{corres2} in the value $\Psi_K^{\circ}(1-\zeta)$ on the norm compatible sequence $1 -\zeta = (1-\zeta_{Np^r})_r$ of $p$-units in the fields $F_r$.  

Finally, we can consider cup products with twisted coefficients.
Let $\omega$ again denote the Teichm\"uller character and $\kappa$ the product of the $p$-adic cyclotomic character with $\omega^{-1}$.
Let $\mc{O}_{Np^r}$
denote the extension of $\zp$ generated by the values of all $\qpbar$-valued characters of $(\Z/Np^r\Z)^{\times}$ for any
$r \ge 1$.
For any even $p$-adic character $\psi$ of $(\Z/Np^s\Z)^{\times}$ (with $s \ge 1$) and $t \in \zp$, we define (as in Section \ref{cyclunits}) 
$$
	\alpha_t^{\psi} = \lim_{r \to \infty}
	\prod_{\substack{i=1\\(i,Np)=1}}^{Np^r-1}
	(1-\zeta_{Np^r}^i)^{\psi\kappa^{t-1}(i)} 
	\in H^1_{\cts}(G_{\Q,S},\mc{O}_{Np^s}(\kappa^t\omega\psi)),
$$
where $\mc{O}_{Np^s}(\kappa^t\omega\psi)$ designates $\mc{O}_{Np^s}$
endowed with a $\kappa^t\omega\psi$-action of $G_{\Q,S}$.  We may then take cup products of pairs of such elements.  Suppose that $k \in \zp$ and that $\theta$ is an
odd character of $(\Z/Np^s\Z)^{\times}$, with the additional assumption that the restriction of $\theta$ to $(\Z/Np\Z)^{\times}$ is primitive.  The cup product 
$\alpha_t^{\psi} \cup \alpha_{k-t}^{\theta\psi^{-1}\omega^{-1}}$ of \eqref{corres3} is then the resulting element of $H^2_{\cts}(G_{\Q,S},\mc{O}_{Np^s}(\kappa^k\omega\theta))$.

On the modular side, we consider the \'etale cohomology group $H^1_{\et}(X_1(Np^r)_{/\qbar};\zp)$.  
Our complex embedding and Poincar\'e duality allow us to identify elements of this Galois module with the singular homology group $H_1(X_1(Np^r);\zp)$ (see Sections \ref{cohomology}-\ref{galact}).  This identifies the $(\pm 1)$-eigenspaces of $H^1_{\et}(X_1(Np^r)_{/\qbar};\zp)$ under complex conjugation with the $(\mp 1)$-eigenspaces of $H_1(X_1(Np^r);\zp)$.  Both of these groups are modules for a cuspidal Hecke algebra, which acts via the adjoint action on cohomology and the standard action on homology,  and we may consider their ordinary parts, i.e., the submodules on which the Hecke operator $U_p$ is invertible.  

The ordinary part of $H_1(X_1(Np^r);\zp)$ contains symbols arising from the classes of paths between cusps in the upper half-plane (see Sections \ref{hom}-\ref{ordinary}).  For $i, j \in \Z$ with $(i,j,Np) = 1$, we may consider the class of the geodesic from $\frac{-b}{dNp^r}$ to $\frac{-a}{cNp^r}$ in homology relative to the cusps, where $ad-bc = 1$, $i \equiv a \bmod Np^r$, and $j \equiv b \bmod Np^r$.  
The symbol $\xi_r(i:j)$ is given by first applying the Manin-Drinfeld splitting to the class of this path and then projecting to the ordinary part.

Inside the part of the cuspidal $\zp$-Hecke algebra that is ordinary and primitive under a certain twisted action of the diamond operators, we have the Eisenstein ideal $I_r$ generated by projections of elements of the form $T_l-1-l\langle l\rangle$ with $l$ prime and $l \nmid Np$, along with $U_l-1$ for $l \mid Np$.  Let $Y_r$ denote the localization of  $H^1_{\et}(X_1(Np^r)_{/\qbar};\zp)$ at the ideal $\mf{m}_r$ generated by $I_r$,
$p$, and $\langle 1+p \rangle-1$, and let $Y_r^-$ denote its $(-1)$-eigenspace under complex conjugation.  In \eqref{corres1}, the symbol $\bar{\xi}_r(i:j)$ then denotes the projection of $\xi_r(i:j)$ to $Y_r^-/I_rY_r^-$ (see Section \ref{cpms}).  

We now define what we shall refer to as two-variable $p$-adic $L$-functions, which are more precisely sequences of Mazur-Tate elements that interpolate such $L$-functions.  
Let $(\Z/Np^r\Z)^{\star}$ denote the set of nonzero elements in $\Z/Np^r\Z$.  If $r \ge 1$ is given, we use $[i]_r$ to denote the element of $\zp[(\Z/Np^r\Z)^{\star}]$ (see Section \ref{modified}) corresponding to $i \in \Z$ with $Np^r \nmid i$.  The $L$-function $\mc{L}_N$ is defined in Section \ref{twovar} as the inverse limit
$$
	\mc{L}_N = 
	\lim_{\substack{\leftarrow\\r}} \sum_{\substack{i=1\\(i,Np)=1}}^{Np^r-1} 
	U_p^{-r}\xi_r(i:1) \otimes [i]_r,
$$
while the modified $L$-function $\mc{L}_N^{\star}$ of Section \ref{modified} is
$$
	\mc{L}_N^{\star} = \lim_{\substack{\leftarrow\\r}} 
	\sum_{i=1}^{Np^r-1} U_p^{-r}\xi_r(i:1) \otimes [i]_r.
$$
The projection of $\mc{L}_N^{\star}$ to the Eisenstein component lies in the completed tensor product $\mc{Y}_N \cotimes \Lambda_N^{\star}$, where $\mc{Y}_N$ denotes the inverse limit of the $Y_r$ and $\Lambda_N^{\star}$ is the inverse limit of the 
$\zp[(\Z/Np^r\Z)^{\star}]$.  The projection of $\mc{L}_N^{\star}$ to 
$\mc{Y}_N^-/\I\mc{Y}_N^- \otimes_{\zp} (\Lambda_N^{\star})^-$ is the object $\overline{\mc{L}_N^{\star}}$  used in \eqref{corres2} (see Section \ref{recLfn}).

We next consider the special values of $\mc{L}_N$.  First, we apply a character of the form $\psi\kappa^{t-1}$, where $t \ge 1$ and $\psi$ is an even character on some $(\Z/Np^s\Z)^{\times}$, obtaining
$$
	\lim_{\substack{\leftarrow\\r}}
	\sum_{\substack{i=1\\(i,Np)=1}}^{Np^r-1} \psi\kappa^{t-1}(i)\xi_r(i:1).
$$
For any odd character $\theta$ on some $(\Z/Np^s\Z)^{\times}$ that is primitive on $(\Z/Np\Z)^{\times}$, we may consider the maximal quotient of the inverse limit of ordinary homology groups 
with $\mc{O}_{Np^s}$-coefficients on which each diamond operator
$\langle j \rangle$ acts as $\theta\omega^{-1}\kappa^{k-2}(j)$.  The image of the above limit in this quotient is denoted $L_p(\xi,\omega\theta,k,\psi,t)$ (see Section \ref{specval}), in that it interpolates the values at the given $t \in \zp$ of the $p$-adic $L$-functions with character $\psi$ of the ordinary cusp forms of weight $k$, level $Np^s$, and character $\theta\omega^{-1}$.
Finally, we may consider its reduction $\overline{L_p(\xi,\omega\theta,k,\psi,t)}$ modulo the Eisenstein ideal of weight $k$ and character $\theta\omega^{-1}$.

The key to relating the above Galois-theoretic and modular objects lies in the construction of maps which take the objects on the left side of our earlier diagram to those on the right side.  These maps should be canonical up to our original choice of complex embedding and make these identifications independent of its choice.  In the paragraph following Proposition \ref{1stisom}, we define, up to a fixed unit in $\Lambda_N$, a homomorphism
$$
	\phi_1 \colon X_K^{\circ} \to \mc{Y}_N^-/\I\mc{Y}_N^-
$$
that arises from the Galois action of $G_{K,S}$ on $\mc{Y}_N$, particularly the map 
$$
	X_K^{\circ} \to \Hom_{\zp}(\mc{Y}_N^+,\mc{Y}_N^-)
$$ 
it induces, together with a modification of a pairing of Ohta's (see Proposition \ref{ohtapair}).
It induces isomorphisms on ``good'' eigenspaces.  
The map yielding \eqref{corres1} is then conjectured to be given by the Tate twist of $\phi_1$ by $\zp(1)$,
and the map yielding \eqref{corres3} is also conjectured to be induced by a twist of $\phi_1$, taking appropriate quotients.

Secondly, we have a homomorphism 
$$
	\phi_2 \colon \mf{X}_K^- \to (\Lambda_N^{\star})^-
$$
determined by the action of $\mf{X}_K^-$ on $p$-power roots of cyclotomic $Np$-units (see Proposition \ref{CKprop}).  
More precisely, $\phi_2(\sigma)$ is the inverse limit of the sequence of elements of $\zp[(\Z/Np^r\Z)^{\star}]$ that have $i$th coefficient modulo $p^s$ given by the exponent of $\zeta_{p^s}$ obtained in applying the Kummer character attached to $\sigma$ to a $p^s$th root of $1-\zeta_{Np^r}^i$.
The map yielding \eqref{corres2} is conjectured to be $\phi_1 \otimes \phi_2$.

\begin{ack}
	The author would like to thank Barry Mazur for his advice and encouragement
	during this project.  He would also like to thank Masami Ohta, William
	Stein, Glenn Stevens, and David Vauclair for doing their best to answer his 
	questions about certain aspects of this work.  
	He also thanks the anonymous referee for a number of helpful comments and suggestions.
	The author's work was supported, in part, by the Canada Research Chairs program
	and by the National Science Foundation under Grant No.\ DMS-0901526.  
	Part of this article was written during stays at the Max Planck Institute for 
	Mathematics, the Institut des Hautes \'Etudes Scientifiques, and the Fields
	Institute, and the author thanks these institutes for their hospitality.
\end{ack}

\section{Galois cohomology} \label{Galois}

\subsection{Iwasawa modules}

\label{p} \label{N}
Let $p$ be an odd prime, and let $N$ be a positive integer prime to $p$.
Let 
\label{F}
$F = \Q(\mu_{Np})$.  Since $\Gal(F/\Q)$ is canonically isomorphic to 
$(\Z/Np\Z)^{\times}$, we may identify characters on the latter group with characters on the former.  Let 
\label{K}
$K$ denote the cyclotomic $\zp$-extension of $F$.  Set 
$$
 \label{ZpN}
	\Z_{p,N} = \lim_{\leftarrow} \Z/N p^r\Z,
$$
and note that that $\Gal(K/\Q)$ is canonically identified 
with $\Z_{p,N}^{\times}$.
Set
$$
 \label{LambdaN}
	\Lambda_N = \zp[[\Z_{p,N}^{\times}]].
$$	
When we speak of $\Lambda_N$-modules, unless stated otherwise,
the action shall be that which arises from the action of $\Gal(K/\Q)$.

We fix, once and for all, a complex embedding 
\label{iota}
$\iota \colon \qbar \hookrightarrow \bf{C}$, which we will use to make a number of canonical choices.
To begin with, for any $d \ge 1$, let 
\label{zetad} 
$\zeta_d = \iota^{-1}(e^{2\pi i/d})$,
which in particular fixes a generator $\zeta = (\zeta_{p^r})$ 
\label{zeta}
of the Tate module.
We use this to identify the Tate module of $K^{\times}$ with $\zp(1)$, though this identification is
primarily notational (e.g., by $\zp(1)$ in a cohomology group, we really mean
the Tate module canonically).

We use 
\label{S}
$S = S_E$ to denote the set of primes dividing $Np$ and any real places in an algebraic extension 
$E$ of $\Q$.  Let  
\label{GES} 
$G_{E,S}$ denote the Galois group of the maximal unramified outside $S$ extension of $E$, and let 
$\mf{X}_E$ denote its maximal abelian pro-$p$ quotient.
\label{OES}
Let $\mc{O}_{E,S}$ denote the ring of $S$-integers of $E$, and let 
 \label{EE} 	
$\mc{E}_E$ denote the $p$-completion of the $S$-units of $E$.
If $\mc{T}$ is a profinite $\zp[[G_{\Q,S}]]$-module and $i \ge 1$, then we let 
$$
\label{HiS} 
	H^i_S(K,\mc{T}) = \lim_{\substack{\leftarrow\\E \subset K}} 
	H^i_{\cts}(G_{E,S},\mc{T}),
$$
in which the inverse limit is taken with respect to corestriction maps
over the number fields $E$ contained in $K$ and containing $F$.  

Let
\label{UK}
$\mc{U}_K$ denote the
group of norm compatible sequences of $S$-units for $K$, i.e., 
$$
	\mc{U}_K = \lim_{\substack{\leftarrow\\ E \subset K}} 
	\mc{O}_{E,S}^{\times} \otimes_{\Z} \zp
	\cong \lim_{\substack{\leftarrow\\ E \subset K}} 
	E^{\times} \otimes_{\Z} \zp.
$$
(Note that any norm compatible sequence must consist of $p$-units, since all decomposition groups
in $\Gal(K/F)$ are infinite and only primes over $p$ ramify, forcing the
valuation of the elements of the sequence to be trivial at primes not over $p$.)
Let
\label{XKS} 
$X_{K,S}$ denote the Galois group of the
maximal abelian pro-$p$ extension of $K$ in which all primes (above those in $S$) split completely.
Kummer theory provides us with the following well-known lemma, of which we sketch a proof for the convenience of the reader.

\begin{lemma}
	There is a canonical isomorphism 
	$$
		H^1_S(K,\zp(1)) \cong \mc{U}_K
	$$
	and a canonical exact sequence 
	\begin{equation} \label{h2seq}
		0 \to X_{K,S} \to H^2_S(K,\zp(1)) \to \bigoplus_{v \in S_K} \zp \to \zp \to 0
	\end{equation}
	of $\Lambda_N$-modules.
\end{lemma}

\begin{proof}
	Let $E$ be an number field in $K$ containing $F$.  Let $A_{E,S}$ denote the
	$p$-part of the $S$-class group of $E$, and let $\mr{Br}_S(E)$ denote the $S$-part
	of the Brauer group of $E$.
	The Kummer sequences arising from the 
	$G_{E,S}$-cohomology of the $S$-unit group of the maximal unramified outside $S$-extension
	of $F$ \cite[Proposition 8.3.11]{nsw} yield
	compatible short exact sequences
	\begin{equation} \label{kummer1}
		0 \to \mc{E}_E/\mc{E}_E^{p^r} \to H^1(G_{E,S},\mu_{p^r}) \to
		A_{E,S}[p^r] \to 0
	\end{equation}
	and
	$$
		0 \to A_{E,S}/p^rA_{E,S} \to H^2(G_{E,S},\mu_{p^r}) \to \mr{Br}_S(E)[p^r] \to 0
	$$
	for $r \ge 1$.
	Considering the isomorphisms
	$$ \label{ctscohom}
		H^i_{\cts}(G_{E,S},\zp(1)) \cong \lim_{\substack{\leftarrow\\r}} 
		H^i(G_{E,S},\mu_{p^r}),
	$$
	the first statement follows from the finiteness of $A_{E,S}$ and the second	by class field theory.
\end{proof}

\subsection{Cup products and the reciprocity map} \label{recmap}
 
Now, consider the cup products
$$
	H^1_{\cts}(G_{E,S},\zp(1)) \otimes_{\zp} 
	H^1_{\cts}(G_{E,S},\zp(1)) \xrightarrow{\cup}
	H^2_{\cts}(G_{E,S},\zp(2))
$$
for number fields $E$ in $K$ containing $F$.  
Note that $\mc{E}_E = \mc{O}_{E,S}^{\times} \otimes_{\Z} \zp$ is canonically isomorphic to
$H^1_{\cts}(G_{E,S},\zp(1))$.  We obtain, therefore, a resulting pairing
$$
\label{cupES}  
	(\,\cdot\,,\,\cdot\,)_{E,S} \colon \mc{E}_E \times \mc{E}_E \to 
	H^2_S(E,\zp(2)).
$$

Recall that $\mc{E}_K$ denotes the $p$-completion of the $S$-units in $K^{\times}$.
In the limit under restriction and corestriction maps, we have
a ``cup product''
$$
	\mc{E}_K \otimes_{\zp} H^1_S(K,\zp(1)) \xrightarrow{\cup} H^2_S(K,\zp(2)),
$$
since $\mc{E}_K$ is canonically isomorphic to the $p$-completion of the
direct limit of the $\mc{E}_E$.
This provides a $\zp$-bilinear pairing
$$
\label{cupKS}
	(\,\cdot\,,\,\cdot\,)_{K,S} \colon \mc{E}_K \times \mc{U}_K \to 
	H^2_S(K,\zp(2)).
$$

\begin{remark}
	In fact, if one takes the limit over $E$ of cup products with 
	$\mu_{p^r}$-coefficients first and then the inverse limit with respect
	to $r$, one obtains a product
	$$
		H^1_{\cts}(G_{K,S},\zp(1)) \otimes_{\zp} H^1_S(K,\zp(1)) 
		\xrightarrow{\cup} H^2_S(K,\zp(2)).
	$$
	The group $H^1_{\cts}(G_{K,S},\zp(1))$ can be identified by Kummer theory
	with the $\Lambda_N$-module with nontrivial elements those
	elements of the $p$-completion 
	of $K^{\times}$ whose $p$-power roots define $\zp$-extensions of $K$ that are 	
	unramified outside $S$.  We shall not need this in this article.
\end{remark}

Consider the exact sequence
\begin{equation} \label{extension}
	1 \to \mf{X}_K \to \mc{T}  \to \zp \to 0
\end{equation}
of $\zp[[G_{K,S}]]$-modules that is determined up to canonical isomorphism by 
the natural projection $\lambda \colon G_{K,S} \to \mf{X}_K$
in the sense that for any lift
$e \in \mc{T}$ of $1 \in \zp$, we have $\lambda(g) = g(e)-e$ for 
all $g \in G_{K,S}$. 
As \eqref{extension} arises as an inverse limit of exact sequences
$$
	1 \to \mf{X}_{F_r} \to \mc{T}_r \to \zp \to 0
$$
given by the projections $\lambda_r \colon 
G_{F_r,S} \to \mf{X}_{F_r}$, we have a 
coboundary map
\begin{equation} \label{coboundary}
	H^1_S(K,\zp) \to H^2_S(K,\mf{X}_K)
\end{equation}
that is the inverse limit of the corresponding coboundaries at the finite level.
For any $r \ge 1$, we have
$$
 	\sigma \cdot \lambda_r(g) = \lambda_r(\sigma g \sigma^{-1})
$$
for $\sigma \in G_{\Q,S}$ and $g \in G_{F_r,S}$, so this is in fact a homomorphism of $\Lambda_N$-modules. 
Twisting \eqref{coboundary} by $\zp(1)$, we obtain a $\Lambda_N$-module
homomorphism 
$$
\label{PsiK} 
	\Psi_K \colon \mc{U}_K \to H^2_S(K,\zp(1)) \otimes_{\zp} \mf{X}_K.
$$
We refer to $\Psi_K$ as the $S$-reciprocity map for $K$.

If $a \in \mc{E}_K$, let 
\label{pia}
$\pi_a \in \Hom_{\cts}(\mf{X}_K,\zp(1))$ denote the corresponding homomorphism.
The cup product relates to $\Psi_K$ as follows:  
\begin{equation} \label{cuppsi}
	 (a,u)_{K,S} = (1 \otimes \pi_a)(\Psi_K(u))
\end{equation}
for $u \in \mc{U}_K$ and $a \in \mc{E}_K$. 

\section{Homology of modular curves} \label{homology}

\subsection{Homology} \label{hom}

\label{p2}
We assume from now on that $p \ge 5$.  
Let $r \ge 1$.
Consider the modular curves 
\label{Y1r}
$Y_1^r(N) = Y_1(Np^r)$ and $X_1^r(N) = X_1(Np^r)$ over ${\bf C}$
and the cusps $C_1^r(N) = X_1^r(N) - Y_1^r(N)$. 
\label{C1r}
We have the following exact sequence in homology:
\begin{equation} \label{homologyseq}
	0 \to H_1(X_1^r(N);\zp) \to H_1(X_1^r(N),C_1^r(N);\zp) \xrightarrow{\delta_r} 
	\tilde{H}_0(C_1^r(N);\zp) \to 0,
\end{equation}
where $\tilde{H}_0$ is used to denote reduced homology.
Let 
\label{Hr}
$\mf{H}_r$ denote the modular Hecke algebra of weight two and level $Np^r$ over $\zp$, 
which acts on $H_1(X_1^r(N),C_1^r(N);\zp)$, and let 
$\mf{h}_r$
\label{hr} 
denote the corresponding cuspidal Hecke algebra over $\zp$, which acts on $H_1(X_1^r(N);\zp)$.  
We have the canonical Manin-Drinfeld splitting over $\qp$,
$$
\label{sr}
	s_r \colon H_1(X_1^r(N),C_1^r(N);\qp) \to H_1(X_1^r(N);\qp).
$$

For any $r \ge 1$, and $a, b \in \Z$ with $(a,b)=1$, let
$$
\label{binom}
	\binom{a}{b}_r \in H_0(C_1^r(N);\zp)
$$
denote the image 
of the cusp corresponding to $a/b \in \bf{P}^1(\Q)$.  
In general, we have that 
$\binom{a}{b}_r = \binom{-a}{-b}_r$ and
\begin{equation} \label{equivcusps}
	\binom{a}{b}_r = \binom{a'+jb'}{b'}_r
\end{equation}
whenever $a \equiv a' \mod Np^r$, $b \equiv b' \mod Np^r$, and $(a,b) = 
(a',b') = 1$ (cf., \cite[Proposition 3.8.3]{ds}).  (We use these equalities
to extend the definition of these symbols to include all $\binom{a}{b}_r$ with
$(a,b,Np)=1$.)

\label{path}
Let $\{\alpha,\beta\}_r$ denote the class in $H_1(X_1^r(N),C_1^r(N);\zp)$
of the geodesic from $\alpha$ to $\beta$ for $\alpha,\beta 
\in \bf{P}^1(\Q)$, which we refer to as a modular symbol.  We note that
the set of such modular symbols 
generate $H_1(X_1^r(N),C_1^r(N);\zp)$ over $\zp$.
They are subject, in particular, to the relations
$$
	\{\alpha,\beta\}_r + \{\beta,\gamma\}_r = \{\alpha,\gamma\}_r
$$
for $\alpha, \beta, \gamma \in \bf{P}^1(\Q)$.
The map $\delta_r$ satisfies
$$
	\delta_r\left(\fpath{a}{c}{b}{d}_r\right) 
	= \binom{b}{d}_r - \binom{a}{c}_r
$$
for $a,b,c,d \in \Z$ with $(a,c) = (b,d) = 1$.

Furthermore, for $u, v \in \Z/Np^r\Z$ with $(u,v) = (1)$, we let
$$
\label{Manin}
	[u:v]_r = \fpath{-b}{dNp^r}{-a}{cNp^r}_r,
$$
where $a, b, c, d \in \Z$ satisfy $ad-bc = 1$,
$u = a \pmod{Np^r}$, and $v = b \pmod{Np^r}$. 
This is the image under the Atkin-Lehner operator
$w_{Np^r}$, which acts
on homology through the matrix 
$$
	\left( \begin{matrix} 0 & -1 \\ Np^r & 0 \end{matrix} \right),
$$ 
of what is usually referred to as a Manin symbol \cite{manin} (i.e., that associated to the pair $(a,b)$).
It is independent of the choices of $a$, $b$, $c$, and $d$.  
We will often abuse notation and refer to $[u:v]_r$ for integers $u$ and $v$
with $(u,v,Np) = 1$.

Recall that $\mf{h}_r$ contains a group of diamond operators 
identified with $(\Z/Np^r\Z)^{\times}$.  We use
\label{<j>r} 
$\langle j \rangle_r$ to denote the element
corresponding to $j \in (\Z/Np^r\Z)^{\times}$.
The homology group $H_1(X_1^r(N),C_1^r(N);\zp)$ 
has a presentation as a $\zp[(\Z/Np^r\Z)^{\times}]$-module with 
generators $[u:v]_r$ for $u, v \in \Z/Np^r\Z$ and $(u,v) = (1)$, subject to the relations:
\begin{gather} 
	[u:v]_r+[-v:u]_r = 0,
	\label{maninrels1}\\
	[u:v]_r = [u:u+v]_r + [u+v:v]_r, 
	\label{maninrels2}\\
	[-u:-v]_r = [u:v]_r, 
	\label{maninrels3}\\
	\langle j \rangle_r^{-1}[u:v]_r = [ju:jv]_r
	\label{maninrels4}
\end{gather}
(see \cite[Theorem 1.9]{manin} for the presentation over $\zp$; the latter
relation is well-known and easily checked).

Additionally, the involution $\alpha \mapsto -\bar{\alpha}$ on the upper half plane provides us with a decomposition of homology into $(\pm 1)$-eigenspaces
\label{Apm}
$H_1(X_1^r(N),C_1^r(N);\zp)^{\pm}$.  We denote the relevant projections of
modular symbols similarly.  
\label{Maninpm}
The presentations of these modules are subject to one additional relation
\begin{equation} \label{maninrels5}
	[-u:v]^{\pm}_r = \pm[u:v]^{\pm}_r.
\end{equation}

\subsection{Ordinary parts} \label{ordinary}

\label{hrord}
The ordinary parts $\mf{h}_r^{\ord}$ and $\mf{H}_r^{\ord}$ of $\mf{h}_r$ and $\mf{H}_r$, respectively, consist of the largest direct summands upon which the $p$th Hecke operator $U_p$ acts invertibly.  Let $e_r$ 
\label{er}
denote Hida's idempotent, which provides maps
\begin{eqnarray*}
	e_r \colon \mf{H}_r \to \mf{H}_r^{\ord} 
	&\mr{and}&
	e_r \colon \mf{h}_r \to \mf{h}_r^{\ord},
\end{eqnarray*}
and similarly for any $\mf{H}_r$-modules.  In particular, \eqref{homologyseq} provides a corresponding exact sequence of ordinary parts:
\begin{equation} \label{ordhomology}
	0 \to H_1(X_1^r(N);\zp)^{\ord} 
	\to H_1(X_1^r(N),C_1^r(N);\zp)^{\ord} 
	\to \tilde{H}_0(C_1^r(N);\zp)^{\ord} \to 0.
\end{equation}
We will identify $H_1(X_1^r(N);\zp)^{\ord}$ with its image in $H_1(X_1^r(N),C_1^r(N);\zp)^{\ord}$.

\begin{lemma} \label{inhomology}
	Suppose that $u, v \in \Z/Np^r\Z$ with $(u,v) = 1$ and 
	both $u$ and $v$ nonzero modulo $p^r$.  Then  
	$$
		e_r[u:v]_r \in H_1(X_1^r(N);\zp)^{\ord}.
	$$
\end{lemma}

\begin{proof}
	By \eqref{ordhomology}, 
	it suffices to show that if $a, b \in \Z$ with $(a,b) = 1$ and
	$p^r \nmid a$, then $e_r\binom{a}{p^rb}_r = 0$.
	This is an immediate corollary of \cite[Proposition 4.3.4]{ohta-ord}.
\end{proof}

For any $u, v \in \Z/Np^r\Z$ with $(u,v) = 1$, let us set
$$
\label{xir}
	\xi_r(u:v) = e_r \circ s_r([u:v]_r).
$$
By Lemma \ref{inhomology}, we have $\xi_r(u:v) = e_r[u:v]_r$ whenever
both $u$ and $v$ are not divisible by $p^r$.

Hida (e.g., \cite{hida-iwa}) 
constructs ordinary Hecke algebras 
\label{h}
\begin{eqnarray*}
	\mf{h} = \lim_{\leftarrow} \mf{h}_r^{\ord} &
	\text{and} & 
	\mf{H} = \lim_{\leftarrow} \mf{H}_r^{\ord}.
\end{eqnarray*}
Inverse limits of the ordinary parts of homology groups with respect to the natural maps of modular curves provide the following $\mf{h}$-modules:
\label{H_1}
\begin{eqnarray*}
	H_1(N) \cong \lim_{\substack{\leftarrow \\ r}} 
	H_1(X_1^r(N);\zp)^{\ord} 
	& \text{and} &
	\mc{H}_1(N) \cong \lim_{\substack{\leftarrow \\ r}} 
	s_r(H_1(X_1^r(N),C_1^r(N);\zp))^{\ord}.
\end{eqnarray*}

We now construct certain inverse limits of our symbols.
\begin{lemma} \label{xidef}
	Let $u \in \Z[\frac{1}{p}]$ and $v \in \Z$.  Suppose that $p \nmid v$ and
	that $(u,v,N)\Z[\frac{1}{p}] = \Z[\frac{1}{p}]$.  Then, for $r$ sufficiently
	large, the symbols $\xi_r(p^r u:v)$ are compatible under
	the natural maps of homology groups, providing an element of 
	$\mc{H}_1(N)$ that we denote 
\label{xi}
	$\xi(u:v)$.
\end{lemma}

\begin{proof}
	Suppose $u = u'p^{-s}$ with $u' \in \Z$ prime to $p$.  Choose $a, b, c, d \in \Z$ with  
	$a \equiv u' \bmod Np^s$, $b \equiv v \bmod Np^r$, and $p^{r-s}ad - bc = 1$.
	Note that 
	$$
		[p^r u: v]_r = \fpath{-b}{dNp^r}{-a}{cNp^s}_r.
	$$ 
	For any $t$ with $s \le t \le r$, this maps to 
	$$
		\fpath{-b}{dNp^r}{-a}{cNp^s}_t = [p^t u :v]_t,
	$$
	since $p^{t-s}a \cdot p^{r-t}d - bc = 1$.
\end{proof}
 
Lemma \ref{inhomology} now has the following immediate corollary.

\begin{corollary}
	Let $u$ and $v$ be as in Lemma \ref{xidef}, and suppose that $u \notin \Z$.
	Then $\xi(u:v) \in H_1(N)$.
\end{corollary}

\subsection{The two-variable $p$-adic $L$-function} \label{twovar}

Mazur \cite{mazur-anom} (but see 
\cite[Section III.2]{mazur-arith}) considers the $\mc{H}_1(N)$-valued measure 
\label{lambdaN}
$\lambda_N$ on $\Z_{p,N}^{\times}$ determined by
$$
	\lambda_N(a + Np^r\Z_{p,N}) = U_p^{-r}\xi(p^{-r}a:1),
$$
where $a \in \Z$ is prime to $Np$ and $r \ge 0$. 
We have an element 
\label{LN}
$\mc{L}_N \in \mc{H}_1(N) \cotimes \Lambda_N$
(where $\cotimes$ denotes the completed tensor product), 
essentially the Mazur-Kitagawa two-variable $p$-adic $L$-function \cite{kitagawa}, determined by
\begin{equation} \label{integrate}
	\tilde{\chi}(\mc{L}_N) = 
	\int_{\Z_{p,N}^{\times}} \chi\lambda_N \in \mc{H}_1(N) \otimes_{\zp} \qpbar
\end{equation}
for any character 
$\chi \in \Hom_{\text{cts}}(\Z_{p,N}^{\times},\qpbar^{\times})$
and induced map 
\label{tildechi}
$$
	\tilde{\chi} \colon \mc{H}_1(N) \cotimes \Lambda_N \to
	\mc{H}_1(N) \otimes_{\zp} \qpbar.
$$
Denoting the group element in $\Lambda_N$ corresponding to 
\label{[j]}
$j \in {\bf Z}_{p,N}^{\times}$ by $[j]$, we have 
$$
	\mc{L}_N = \lim_{\leftarrow} \sum_{\substack{j=0\\(j,Np) = 1}}^{Np^r-1} 
	U_p^{-r}\xi_r(j:1)
	\otimes [j]_r \in \mc{H}_1(N) \cotimes \Lambda_N,
$$
where $[j]_r$ 
\label{[j]r}
denotes the image of $[j]$ in $\zp[(\Z/Np^r\Z)^{\times}]$.

We shall require certain modified versions of this $L$-function.  In this section, we mention the following generalization. For any $M$ dividing $N$, let us set
\begin{equation} \label{LNM}
	\mc{L}_{N,M} = 
	\lim_{\leftarrow} \sum_{\substack{j=0\\(j,Np) = 1}}^{Np^r-1} 
	U_p^{-r}\xi_r(j:M)
	\otimes [j]_r \in \mc{H}_1(N) \cotimes \Lambda_N.
\end{equation}
One can also define this similarly to \eqref{integrate} by integration, replacing $\lambda_N$ by 
\label{lambdaNM}
$\lambda_{N,M}$ with
$$
	\lambda_{N,M}(a + Np^r\Z_{p,N}) = U_p^{-r}\xi(p^{-r}a:M).
$$
We will now explain why $\mc{L}_{N,M}$ is well-defined. 
 
In general, suppose that $t$ is a positive divisor of $Np^r$ for some $r$ and $u$ and $v$ are positive integers not divisible by $Np^r$ with $(tu,v,Np) = 1$. Let $Q = Np^r/t$, and choose $a,b,c,d \in \Z$ with $tad-bc = 1$, $a \equiv u \bmod Np^r$,
and $b \equiv v \bmod Np^r$.  For any such $t$, we define
\label{Ut}
$$
	U_t = \prod_{\substack{l \mid Np\\l \text{ prime}}} U_l^{m_l},
$$
where $m_l$ denotes the $l$-adic valuation of $t$.  Then
\begin{align*}
	U_t[tu:v]_r &= U_t
		\fpath{-b}{dNp^r}{-ta}{cNp^r}_r
	\\&= \sum_{k=0}^{t-1}
		\fpath{-b+kdNp^r}{tdNp^r}{-a+kcQ}{cNp^r}_r	
	\\
	&= \sum_{k=0}^{t-1}[u+kQ:v]_r.
\end{align*}
We obtain
\begin{equation} \label{maninsum}
	U_t\xi_r(tu:v) = \sum_{k=0}^{t-1} \xi_r(u+kQ:v).
\end{equation}

Hence, for $s \ge r$ and
any positive $i < Np^r$ with $(i,Np)=1$, the quantity
$$
	\sum_{k=0}^{p^{s-r}-1} \xi_s(i+kNp^r:M) \otimes [i+kNp^r]_s
	\in H_1(X_1^s(N),C_1^s(N);\zp)^{\ord} \otimes_{\zp} 
	\zp[(\Z/Np^s\Z)^{\times}]
$$
maps to
$$
	U_p^{s-r}\xi_s(p^{s-r}i:M) \otimes [i]_r \in 
	H_1(X_1^s(N),C_1^s(N);\zp)^{\ord} \otimes_{\zp} 
	\zp[(\Z/Np^r\Z)^{\times}]
$$ 
and, therefore, to
$$
	U_p^{s-r}\xi_r(i:M) \otimes [i]_r \in 
	H_1(X_1^r(N),C_1^r(N);\zp)^{\ord} \otimes_{\zp} 
	\zp[(\Z/Np^r\Z)^{\times}]
$$
under the maps inducing the inverse limit in \eqref{LNM}.

The following is immediate from Lemma \ref{inhomology}.

\begin{corollary}
	The $L$-function $\mc{L}_{N,M}$ lies in $H_1(N) \cotimes \Lambda_N$.
\end{corollary}

\subsection{Cohomology} \label{cohomology}

We now explore the relationship between homology and cohomology groups
of modular curves.  We show that, for our purposes, they are interchangeable.
For this, we consider exact sequences in reduced (singular) homology and cohomology of our modular curves and commutative diagrams induced by Poincar\'e
duality.  We refer the reader to \cite[Section 1.8]{stevens}
as well.

\begin{proposition} \label{poincare}
	For $r \ge 1$, we have canonical commutative diagrams
	\begin{equation} \label{lefschetz}
		\SelectTips{cm}{} \xymatrix@C=15pt{
		0 \ar[r] & H_1(X_1^r(N);\zp) \ar[r] \ar[d]^{\wr} & 
		H_1(X_1^r(N),C_1^r(N);\zp) \ar[r] \ar[d]^{\wr} &
		\tilde{H}_0(C_1^r(N);\zp) \ar[r] \ar[d]^{\wr} & 0\\
		0 \ar[r] & H^1(X_1^r(N);\zp) \ar[r] &
		H^1(Y_1^r(N);\zp) \ar[r] &
		\tilde{H}^0(C_1^r(N);\zp) \ar[r] & 0
		}
	\end{equation}
	that are compatible with the natural maps on homology and trace maps on
	cohomology.  Furthermore, the actions of $\mf{H}_r$ 
	on the homology groups and the adjoint Hecke algebras
\label{Hr*}
	$\mf{H}_r^*$ on the cohomology groups are compatible.  
\end{proposition}

\begin{proof}
	Let $D = \zp[{\bf P}^1(\Q)]$, and let $D_0$ denote the kernel of the obvious
	augmentation map $D \to \zp$.  Set $G_r = \Gamma_1(Np^r)$.
	Using the homological version of \cite[Proposition 4.2]{as},
	we may rewrite the top exact sequence in \eqref{lefschetz} canonically
	as
	$$
		0 \to \ker \alpha  \to (D_0)_{G_r} \xrightarrow{\alpha} 
		\ker (D_{G_r} \to \zp) \to 0.
	$$
	As in \cite[loc.\ cit.]{as}, the $\zp$-dual of
	this sequence is canonically
	\begin{equation} \label{dualseq}
		0 \leftarrow H^1(X_1^r(N);\zp) \leftarrow H^1_c(Y_1^r(N);\zp)
		\leftarrow \tilde{H}^0(C_1^r(N);\zp) \leftarrow 0
	\end{equation}
	as an exact sequence of $\mf{H}_r$-modules.  Finally, Poincar\'e 
	duality implies that the $\zp$-dual of the latter sequence is
	canonically the exact sequence of $\mf{H}_r^*$-modules,
	$$
		0 \to H^1(X_1^r(N);\zp) \to H^1(Y_1^r(N);\zp) \to 
		\tilde{H}^0(C_1^r(N);\zp) \to 0,
	$$
	via cup product (fixing a generator of $H^2_c(Y_1^r(N);\zp)$
	corresponding to a simple counterclockwise loop around a point in the 
	upper half-plane), and it is well-known that the Hecke and
	adjoint Hecke actions are compatible with the cup product.
	
	Now, the natural surjections $(D_0)_{G_s} \to (D_0)_{G_r}$ for $s \ge r$ 
	yield the natural injections
	$$
		\Hom_{G_r}(D_0,\zp) \to \Hom_{G_s}(D_0,\zp),
	$$ 
	in the dual, and these maps are all compatible with the standard 
	Hecke actions arising from the action of $GL_2(\Q)^+$ on $D_0$. 
	Furthermore, the
	trace maps $H^1(Y_1^s(N);\zp) \to H^1(Y_1^r(N);\zp)$ are compatible
	with the actions of the adjoint Hecke algebras, and agree with 
	the latter inclusions under Poincar\'e duality.	
	The rest follows easily.
\end{proof}

As before, we have a Manin-Drinfeld splitting
$$
	s^r \colon H^1(Y_1^r(N);\qp) \to H^1(X_1^r(N);\qp)
$$
and cohomology groups
\label{H^1}
\begin{eqnarray*}
	H^1(N) \cong \lim_{\substack{\leftarrow \\ r}} 
	H^1(X_1^r(N);\zp)^{\ord} 
	& \text{and} &
	\mc{H}^1(N) \cong \lim_{\substack{\leftarrow \\ r}} 
	s^r(H^1(Y_1^r(N);\zp))^{\ord},
\end{eqnarray*}
where the inverse limits are taken with respect to trace maps and
\label{ordcohomology}
``ord'' now denotes the part upon which 
	the adjoint Hecke operator
$U_p^*$ acts invertibly.
These are modules over
\label{h*}
\begin{eqnarray*}
	\mf{h}^* = \lim_{\leftarrow} (\mf{h}_r^*)^{\ord} &
	\text{and} & 
	\mf{H}^* = \lim_{\leftarrow} (\mf{H}_r^*)^{\ord},
\end{eqnarray*}
respectively.

\subsection{Galois actions} \label{galact}

Our fixed embedding $\iota \colon \qbar \hookrightarrow {\bf C}$ defines compatible isomorphisms
$$
	\Phi^r \colon H^1(X_1^r(N);\qp) 
	\xrightarrow{\sim} H^1_{\et}(X_1^r(N){}_{/\qbar};\qp)
$$
and therefore an isomorphism $\Phi$ in the inverse limit.
We define
\label{H^1et}
\begin{eqnarray*}
	H^1_{\et}(N) = \Phi(H^1(N)) &\mr{and}& \mc{H}^1_{\et}(N)
	= \Phi(\mc{H}^1(N)).
\end{eqnarray*}

Using, for instance, the duality between the top sequence in \eqref{lefschetz} and the exact sequence in \eqref{dualseq}, we have Galois actions on homology as well, producing \'etale homology groups and isomorphisms
$$
	\Phi_r \colon H_1(X_1^r(N);\qp) \xrightarrow{\sim}
	H_1^{\et}(X_1^r(N){}_{/\qbar};\qp),
$$
resulting in an isomorphism in the inverse limit that we also label $\Phi$.
We define
\label{H_1et}
\begin{eqnarray*}
	H_1^{\et}(N) = \Phi(H_1(N)) &\mr{and}& \mc{H}_1^{\et}(N)
	= \Phi(\mc{H}_1(N)).
\end{eqnarray*}

Note that the isomorphisms between \'etale homology and cohomology groups resulting from Proposition \ref{poincare} and our choice of $\iota$ are not isomorphisms of Galois modules.  
Rather, Poincar\'e duality yields a perfect pairing
$$
	H^1_{\et}(X_1^r(N){}_{/\qbar};\zp) \times
	H^1_{\et}(X_1^r(N){}_{/\qbar};\zp(1)) \to 
	H^2_{\et}(X_1^r(N){}_{/\qbar};\zp(1)) \cong \zp,
$$
of Galois modules.  
We have canonical isomorphisms $H_1^{\et}(N) \cong H^1_{\et}(N)(1)$, and similarly, $\mc{H}_1^{\et}(N) \cong \mc{H}^1_{\et}(N)(1)$.  Though we will continue identify elements of $\mc{H}_1^{\et}(N)$ with elements of $\mc{H}^1_{\et}(N)$, we also need to remain aware of the Galois actions for later applications.

Note that the image of $\mc{L}_{N,M}$ in $\mc{H}_1^{\et}(N) \cotimes \Lambda_N$ depends upon $\iota$, since $\Phi_r$ applied to $\xi_r(j:M)$ for $j$ prime to $Np$ varies with $\iota$ (i.e., is not fixed by the absolute Galois group $G_{\Q}$).  

\section{First form of the conjecture} \label{conjecture}

\subsection{Eigenspaces}

\label{p3}
We continue to fix $p$ prime (with $p \ge 5$) and $N \ge 1$ prime to $p$. 
\label{N2}
 We assume from now on that $(\Z/N\Z)^{\times}$ has prime-to-$p$ order.  That is, we assume that $p$ does not divide $\varphi(N)$, where $\varphi$ denotes the Euler-phi function.

For a $\zp[(\Z/Np\Z)^{\times}]$-module $A$, we define the primitive part of $A$ to be
$$
	\ker\biggl( A \to \bigoplus_{\substack{M \mid Np\\
	Np/M\, \text{prime}}} A
	\otimes_{\zp[(\Z/Np\Z)^{\times}]} \zp[(\Z/M\Z)^{\times}] \biggr),
$$
Since $p \nmid \varphi(N)$, the primitive part of $A$ is canonically a
direct summand of $A$ with complement 
$$
	\sum_{\substack{M \mid Np\\
	Np/M\, \text{prime}}} A^{\ker((\Z/Np\Z)^{\times} \to (\Z/M\Z)^{\times})}.
$$
We define 
\label{Acirc}
$A^{\circ}$ to be the submodule of $A$ consisting of all elements of the primitive part of $A$ upon which $-1 \in (\Z/Np\Z)^{\times}$ acts as multiplication by $-1$, i.e, the odd part of the primitive part of $A$.

\begin{remark}
	For now, we work with the above definition of $A^{\circ}$.  Later, 
	the notation $A^{\circ}$ will depend upon $A$.
\end{remark}

We may phrase this in terms of eigenspaces of $\zp[(\Z/Np\Z)^{\times}]$-modules.  
Given a Dirichlet character $\chi \colon (\Z/Np\Z)^{\times} \to \qpbar^{\times}$ of conductor dividing $Np$, let 
\label{Rchi}
$R_{\chi}$ denote the ring generated over $\zp$ by the values of $\chi$.
Then $R_{\chi}$ is canonically a quotient of
$\zp[(\Z/Np\Z)^{\times}]$. For a $\zp[(\Z/Np\Z)^{\times}]$-module $A$, set 
$$
\label{A(chi)}
	A^{(\chi)} = A \otimes_{\zp[(\Z/Np\Z)^{\times}]} R_{\chi}.
$$ 
This is canonically a quotient of $A$ and is an $R_{\chi}[(\Z/Np\Z)^{\times}]$-module with a $\chi$-action.  

Let 
\label{Sigma}
$\Sigma$ denote the set of $G_{\qp}$-conjugacy classes of Dirichlet characters on $(\Z/Np\Z)^{\times}$.  We use $(\chi)$ to denote the class of
$\chi$.  The direct sum of the quotient maps gives rise to a decomposition
$$
	A \cong \bigoplus_{(\chi) \in \Sigma} A^{(\chi)},
$$
canonical up to the choice of representatives of the classes.  Let $\Sigma_{Np}$ denote the subset of $\Sigma$ consisting of classes of primitive characters, i.e., of characters of conductor $Np$.
\label{SigmaNp}
For a $\zp[(\Z/Np\Z)^{\times}]$-module $A$, we then have
$$
	A^{\circ} \cong \bigoplus_{\substack{(\chi) \in \Sigma_{Np}\\ 
	\chi \text{ odd}}} A^{(\chi)}.
$$

\subsection{Eisenstein components}

We will have need to distinguish between Galois and Hecke actions of $\Lambda_N$ on certain modules
that have both.  Therefore, we write
\label{LambdaNh}
$\Lambda_N^{\mf{h}}$ for $\Lambda_N$ when we consider it together with its canonical surjection onto 
the $\zp$-subalgebra of $\mf{h}$
(resp., $\mf{h}^*$) topologically generated by the diamond operators \label{<j>}
$\langle j \rangle$ (resp., adjoint diamond operators $\langle j \rangle^*$) 
\label{<j>*} for $j \in \Z_{p,N}^{\times}$. 
Let
\label{varepsilon}
$\varepsilon \colon \Lambda_N^{\mf{h}} \to (\Lambda_N^{\mf{h}})^{\circ}$ denote the natural projection map, 
viewing $\Lambda_N^{\mf{h}}$ as a $\zp[(\Z/Np\Z)^{\times}]$-module in the obvious manner.
Let 
\label{omega}
$\omega \colon (\Z/Np\Z)^{\times} \to \zp^{\times}$ denote the Dirichlet (Teichm\"uller) character which factors 
as projection to $(\Z/p\Z)^{\times}$ followed by the natural inclusion,
and which we will also view as a character on $\Z_{p,N}^{\times}$.
Let 
\label{kappa}
$\kappa \colon \Z_{p,N}^{\times} \to \zp^{\times}$ denote the canonical projection to $1+p\zp$.
We define the Eisenstein ideal 
\label{I}
$\I$ of $\mf{h}$ to be the ideal generated by 
$T_l - 1 - l\langle l \rangle$  and $\langle l \rangle - \varepsilon(\langle l \rangle)\omega(l)^{-1}$
for $l \nmid Np$, along with $U_l - 1$ for $l \mid Np$.  
Let $\mf{m} = \I + (p,\langle 1+p \rangle -1)\mf{h}$.
(Despite the notation, $\mf{m}$ need not be a maximal ideal of $\mf{h}$.)
Using the same definition with adjoint operators, we have corresponding ideals of $\mf{h}^*$,
which we also denote $\I$ and $\mf{m}$, respectively, by abuse of notation.
We also have an Eisenstein ideal $\mf{I}$ of $\mf{H}$ with the same generators and 
$\mf{M} = \mf{I} + (p,\langle 1+p \rangle -1)\mf{H}$
\label{M}
(and similarly for $\mf{H}^*$).

We define the ``localization'' of $\mf{h}^*$ at $\mf{m}$ by
$$
\label{hm}
	\mf{h}^*_{\mf{m}} = \prod_{\substack{\mf{m}' \subset \mf{h}^*	
	\text{ maximal}\\ \mf{m} 	\subset \mf{m}'}} \mf{h}^*_{\mf{m'}}
$$ 
(and similarly for $\mf{h}$).
This is well-known to be a direct summand of $\mf{h}^*$.  
When we refer to elements of $\mf{h}^*$ and its modules as elements of $\mf{h}^*_{\mf{m}}$ and localizations at $\mf{m}$ of said modules, we shall mean after taking the appropriate projection map.  We use this notation without further comment. 

When needed, we will denote the eigenspace of an $\h$-module $Z$ upon which the 
adjoint diamond operators $\langle j \rangle^*$ with $j \in
(\Z/Np\Z)^{\times}$ act by $\theta\omega^{-1}(j)$, where $\theta$ is a primitive, odd Dirichlet character of conductor $Np$, by 
\label{Ztheta}
$Z^{\ltheta}$.  We remark that $\mf{m}^{\ltheta}$ is the 
maximal ideal of 
the nontrivial eigenspace $(\h)^{\ltheta}$ when $p$ divides the 
generalized Bernoulli number 
\label{B1theta}
$B_{1,\theta}$, 
and the inverse images of such $\mf{m}^{\ltheta}$ in $\mf{h}^*$ are exactly 
the maximal ideals of $\mf{h}^*$ containing $\mf{m}$.
	
Let 
\label{ZN}
$\mc{Z}_N = \mc{H}^1_{\et}(N)_{\mf{m}}$ and 
$\mc{Y}_N = H^1_{\et}(N)_{\mf{m}}$.  
	We note the following useful fact.

	\begin{lemma} \label{usefulfact}
	The inverse limit of the maps $s^r$ induces a canonical isomorphism
	$$
		(\lim_{\leftarrow} H^1_{\mr{\acute{e}t}}(Y_1^r(N);\zp)^{\ord})
		\otimes_{\mf{H}^*} \h \xrightarrow{\sim} \mc{Z}_N.
	$$
	\end{lemma}

	\begin{proof}
		The action of $\mf{H}^*_{\mf{M}}$ on $\mc{Z}_N$ factors by definition through
		the action of $\mf{h}^*_{\mf{m}}$, so the $s^r$ do indeed induce a canonical surjective map 
		$s$ as in the statement of the lemma, which we must show is injective.  
		For this, we first note that the natural map $\iota \colon \mc{Y}_N \to \mc{Z}_N$ is by 
		definition injective, as is then the natural map 
		$$
			t \colon \mc{Y}_N \to (\lim_{\leftarrow} H^1_{\et}(Y_1^r(N);\zp)^{\mr{ord}})
			\otimes_{\mf{H}^*} \h,
		$$ 
		given that $s \circ t = \iota$.  By \cite[Theorem 1.5.5]{ohta-cong} and \cite[Corollary
		A.2.4]{ohta-cong}, we have that
		the congruence module $\mc{Z}_N/\mc{Y}_N$ is isomorphic to $\h/\I$.  In turn, this is
		isomorphic by \cite[Theorem 2.3.6]{ohta-cong} to
		$$
			\mf{H}^*_{\mf{M}}/\mf{I} \otimes_{\mf{H}^*} \h \cong 
			(\lim_{\leftarrow} \tilde{H}^0_{\et}(C_1^r(N),\zp)^{\mr{ord}}) \otimes_{\mf{H}^*} \h,
		$$
		which is canonically the cokernel of $t$.  It follows that $s$ is injective as well.
	\end{proof}

We have that $\mc{Z}_N$ (resp., $\mc{Y}_N$) decomposes into a direct sum of $(\pm 1)$-eigenspaces 
\label{Apm2}
$\mc{Z}_N^{\pm}$ (resp., $\mc{Y}_N^{\pm}$) under the complex conjugation determined by our complex
embedding $\iota$.  We wish to compare this decomposition to another standard sort of decomposition, determined locally at a prime above $p$, that is well-understood by work of Ohta \cite{ohta-comp},
building on work of Mazur-Wiles \cite{mw2} and Tilouine \cite{tilouine}.

Let $\theta$ be a primitive, odd Dirichlet character of conductor $Np$ with $p \mid B_{1,\theta}$.
For now, let
$D_p$ be an arbitrary decomposition group at $p$, and let 
$\delta$ be any element of its inertia subgroup $I_p$
such that $\omega(\delta)$ has order $p-1$ and such that the closed subgroup 
generated by $\delta$
has trivial maximal pro-$p$ quotient.  
Set 
\label{Zpt}
$Z_{\theta} = 
(\mc{Z}_N^{\ltheta})^{I_p}$, and let 
$Z'_{\theta}$ consist of those elements of $\mc{Z}_N^{\ltheta}$ upon which $\delta$ acts
by a nontrivial power of $\omega(\delta)$.
\label{Zmt}

Abusing notation,
we denote 
the eigenspace of $\Lambda_N^{\mf{h}}$ upon which the group element $[j]$ for $j \in (\Z/Np\Z)^{\times}$
acts by $\theta\omega^{-1}$ by
\label{LambdaNt}
$\Lambda_N^{\ltheta}$, and we use 
\label{LNt}
$\mf{L}_N^{\ltheta}$ to denote its quotient field.  Recalling, for instance, 
\cite[Corollary 2.3.6]{ohta-ord2} 
and \cite[Lemma 5.1.3]{ohta-ord} and using the fact that $(\mc{Z}_N/\mc{Y}_N)^{\ltheta}$ is 
$\Lambda_N^{\ltheta}$-torsion, 
we have that $\mc{Z}_N^{\ltheta} = Z'_{\theta} \oplus Z_{\theta}$ as $(\h)^{\ltheta}$-modules, where
$Z_{\theta}$ is free of rank $1$ and the tensor product of $Z'_{\theta}$ with
$\mf{L}_N^{\ltheta}$ over $\Lambda_N^{\ltheta}$ is free over 
$(\h)^{\ltheta} \otimes_{\Lambda_N} \mf{L}_N^{\ltheta}$.

Consider the representation
\label{rhot}
$$
	\rho_{\theta} \colon G_{\Q} \to \mr{Aut}_{\h}(\mc{Z}_N^{\ltheta})
$$
of the absolute Galois group $G_{\Q}$, and the four maps 
\begin{eqnarray*} 
\label{dt} 
	a_{\theta} \colon G_{\Q} \to \End_{\h}(Z'_{\theta}), &&
	b_{\theta} \colon G_{\Q} \to \Hom_{\h}(Z_{\theta},Z'_{\theta}),\\ 
	c_{\theta} \colon G_{\Q} \to \Hom_{\h}(Z'_{\theta},Z_{\theta}), &&
	d_{\theta} \colon G_{\Q} \to \End_{\h}(Z_{\theta})
\end{eqnarray*}
that $\rho_{\theta}$ induces, which allow us to view $\rho_{\theta}$ in matrix form as
$$
		\rho_{\theta}(\sigma) = \left( \begin{matrix} a_{\theta}(\sigma) 
		&b_{\theta}(\sigma) \\ c_{\theta}(\sigma) & d_{\theta}(\sigma) 
		\end{matrix} \right)
$$
for $\sigma \in G_{\Q}$.
Note that $\End_{\h}(Z_{\theta}) \cong (\h)^{\ltheta}$ and similarly for $Z'_{\theta}$, and let
\label{Bt} 
$B_{\theta}$ and $C_{\theta}$ denote the $(\h)^{\ltheta}$-modules generated by
the images of $b_{\theta}$ and $c_{\theta}$
respectively.
	
\begin{proposition} \label{rhoimage}
	Let $\theta$ be a primitive, odd, non-quadratic Dirichlet character of conductor $Np$
	such that $p \mid B_{1,\theta}$.  Then 
	\begin{equation} \label{image}
		\rho_{\theta}(G_K) = 
		 \left\{ \left( \begin{matrix} \alpha & \beta \\ \gamma & \delta 
		\end{matrix} \right) \mid
		\alpha, \delta \in 1 + \I^{\ltheta},\, \beta \in B_{\theta},\, \gamma \in C_{\theta}, \,
		\alpha\delta-\beta\gamma = 1 
		\right\}.
	\end{equation}
\end{proposition}
	
\begin{proof}
	The 
	induced maps $\bar{b}_{\theta} \colon G_K \to B_{\theta}/\I B_{\theta}$
	and $\bar{c}_{\theta} \colon G_K \to C_{\theta}/\I C_{\theta}$ 
	are surjective homomorphisms, which follows as in
	\cite[Lemma 5.3.18]{ohta-ord} 
	(with $\mc{Z}_N$ replacing $\mc{Y}_N$ and for $C_{\theta}$ just as for $B_{\theta}$).  
	Since $\theta^2 \neq 1$,
	eigenspace considerations yield that
	the fixed fields of the kernels of $\bar{b}_{\theta}$ and $\bar{c}_{\theta}$ on $G_K$
	intersect precisely in $K$.
	
	Let $G$ denote the group on the right-hand side of \eqref{image}, which we know
	contains $\rho_{\theta}(G_K)$ by the same argument as in \cite[Lemma 5.3.12]{ohta-ord}
	(as described for instance in \cite[Section 4.2]{ohta-comp2}, noting Lemma \ref{usefulfact}).
	To prove the proposition, it suffices to note that the diagram 
	$$
		\SelectTips{cm}{}  \xymatrix{
		G_K \ar[r]^{\rho_{\theta}} \ar@{>>} [rd]_-{\bar{b}_{\theta} + \bar{c}_{\theta}} 
		& G \ar@{>>} [d]^-{\left( \begin{smallmatrix} \alpha & \beta \\ \gamma & \delta \end{smallmatrix} \right) 		\mapsto \bar{\beta} + \bar{\gamma}} \\
		& B_{\theta}/\I B_{\theta} \oplus C_{\theta}/\I C_{\theta}
		}
	$$
	commutes, with the vertical map inducing an isomorphism on $G^{\mr{ab}}$.
	For this latter claim, we compute the commutator subgroup of $G$.  
	
	Note that $B_{\theta}C_{\theta} = \I^{\ltheta}$ by the same argument 
	leading to \cite[Corollary 5.3.13]{ohta-ord}, 
	together with \cite[Corollary 4.1.12]{ohta-comp2}.
	Let us set 
	$$
		G' =  \left\{ \left( \begin{matrix} \alpha & \beta \\ \gamma & \delta 
		\end{matrix} \right) \mid
		\alpha, \delta \in 1 + \I^{\ltheta},\, 
		\beta \in \I B_{\theta},\, \gamma \in \I C_{\theta}, \,
		\alpha\delta-\beta\gamma = 1 
		\right\},
	$$
	and let $L$, $D$, and $U$ denote the subgroups of $G'$ consisting of lower-triangular
	unipotent, diagonal, and upper-triangular unipotent matrices in $G'$, respectively.
	Since $G' = LDU$ as sets, our claim amounts to showing that $L$, $D$, and $U$ 
	are subgroups of $[G,G]$.
	
	First, note that
	$$
		\left[ \left(  \begin{matrix} \alpha & 0 \\ 0 & \alpha^{-1}
		\end{matrix} \right) , \left(  \begin{matrix} 1 & \beta \\ 0 & 1
		\end{matrix} \right) \right] 
		= \left(  \begin{matrix} 1 & (\alpha^2-1)\beta \\ 0 & 1
		\end{matrix} \right) 
	$$
	for $\alpha \in 1+\I^{\ltheta}$ and $\beta \in B_{\theta}$,
	so $U \subset [G,G]$.
	Similarly, we have $L \subset [G,G]$.  
	Next, let $\beta \in B_{\theta}$ and $\gamma \in C_{\theta}$, 
	set $t = \beta\gamma \in \I^{\ltheta}$, and
	consider the commutator
	$$
		\left[  \left(  \begin{matrix} 1 & 0 \\ \gamma & 1
		\end{matrix} \right), 
		\left(  \begin{matrix} 1&  \beta \\ 0 & 1
		\end{matrix} \right) \right] 
		= \left(  \begin{matrix} 1-t & t\beta \\ 
		-t\gamma & 1+t+t^2
		\end{matrix} \right).
	$$
	Setting $u = 1-t$, and observing that
	$$
		\left(  \begin{matrix} u &  0 \\ 0 & u^{-1}
		\end{matrix} \right) =  \left(  \begin{matrix} 1 & 0 \\ tu^{-1}\gamma & 1
		\end{matrix} \right)
		\left(  \begin{matrix} 1-t& t\beta \\ 
		-t\gamma & 1+t+t^2
		\end{matrix} \right)
		\left(  \begin{matrix} 1& - tu^{-1}\beta \\ 0 & 1
		\end{matrix} \right),
	$$
	we have $D \subseteq [G,G]$.  Thus $G' = [G,G]$, and $G^{\mr{ab}}$ is as desired.
\end{proof}

We are now ready to compare the two types of decompositions of $\mc{Z}_N^{\ltheta}$.

\begin{theorem} \label{decomp}
	For any primitive, odd Dirichlet character $\theta$ of conductor $Np$ such that $p \mid B_{1,\theta}$,
	there exist a decomposition group $D_p$ at $p$ in $G_{\Q}$ and an element $\delta$ of
	its inertia subgroup $I_p$ for which $\omega(\delta)$ has order $p-1$ and the closed subgroup 
	generated by $\delta$ has trivial pro-$p$ quotient such that 
	$(\mc{Z}_N^+)^{\ltheta} = (\mc{Z}_N^{\ltheta})^{I_p}$ and $(\mc{Z}_N^-)^{\ltheta}$ is the
	submodule of $\mc{Z}_N^{\ltheta}$ upon which $\delta$ acts by a nontrivial power of $\omega(\delta)$. 
 \end{theorem}

\begin{proof}
	It suffices to show that for any choice of $D_p$ and $\delta$, there exists a conjugate 
	$\tau$ of our fixed complex conjugation $\tau_0$ such that, in our earlier notation,
	$\rho_{\theta}(\tau) = \left( \begin{smallmatrix} -1 & 0 \\ 0 & 1 \end{smallmatrix} \right)$.
	Letting $\overline{\det \rho_{\theta}}$ denote the composition of $\det \rho_{\theta}$
	with the projection $(\h)^{\ltheta} \to (\mf{h}^*/\mf{m})^{\ltheta}$, the image of $\rho_{\theta}$
	is isomorphic to the semi-direct product of the pro-$p$ group 
	$\ker(\overline{\det \rho_{\theta}})/\ker \rho_{\theta}$ with the finite prime-to-$p$ group 
	$\mr{im}(\overline{\det \rho_{\theta}})$ (see, for example, \cite[Lemma 3.3.5]{ohta-ord2}).  
	
	Since any two Sylow $2$-subgroups in the image of $\rho_{\theta}$ 
	are conjugate, two elements of order two in the
	image are conjugate by the image of an element of $G_{\Q}$ if their
	determinants agree.  As $\det \rho_{\theta}(\tau_0) = -1$, it suffices to show that 
	$\left( \begin{smallmatrix} -1& 0 \\ 0& 1 \end{smallmatrix} \right) \in \rho_{\theta}(G_{\Q})$.
	If $\theta(\delta)$ has even order $2m$, then
	$\rho_{\theta}(\delta^m) = \left( \begin{smallmatrix} -1& 0 \\ 0& 1 \end{smallmatrix} \right)$.	
		
	If $\theta(\delta)$ has odd order, which in particular implies $\theta^2 \neq 1$,
	we take a different approach, exploiting our knowledge of the image of $\rho_{\theta}$.
	Taking tensor products over $\Lambda = \Lambda_N^{\ltheta}$, we have
	$$
		(\mc{Z}_N^{\pm})^{\ltheta} \otimes_{\Lambda} \mf{L}_N^{\ltheta} 
		\cong Z_{\theta} \otimes_{\Lambda} \mf{L}_N^{\ltheta}
		\cong Z'_{\theta} \otimes_{\Lambda} \mf{L}_N^{\ltheta}
		\cong (\h)^{\ltheta} \otimes_{\Lambda} \mf{L}_N^{\ltheta}
	$$
	by \cite[Lemma 5.1.3]{ohta-ord} and \cite[p.\ 588]{hida}.
	We therefore have that
	$\rho_{\theta}(\tau_0) =
	 P\left( \begin{smallmatrix} -1& 0 \\ 0& 1 \end{smallmatrix} \right)P^{-1}$ for
	some
	$$
		P \in \mr{Aut}_{\h}(\mc{Z}_N^{\ltheta} \otimes_{\Lambda} \mf{L}_N^{\ltheta}) \cong
		GL_2((\h)^{\ltheta} \otimes_{\Lambda} \mf{L}_N^{\ltheta})
	$$ 
	of determinant $1$.
	If
	$P =  \left( \begin{smallmatrix} \alpha & \beta \\ \gamma & \delta 
		\end{smallmatrix} \right)$,
	then 
	$$
		\rho_{\theta}(\tau_0) 
		=  \left( \begin{matrix} -(\alpha\delta+\beta\gamma) & 2\alpha\beta 
		\\ -2\gamma\delta & \alpha\delta+\beta\gamma \end{matrix} \right).
	$$
	Since $\alpha\delta-\beta\gamma = 1$ and $\alpha\beta\gamma\delta \in \I^{\ltheta}$, exactly one of
	$\beta\gamma$ and $\alpha\delta$ is a unit in $(\h)^{\ltheta}$.  
	However, it cannot be $\beta\gamma$,
	as this would force 
	$$
		\alpha\delta+\beta\gamma \equiv \beta\gamma \equiv 1 \bmod \mf{m}^{\ltheta},
	$$
	contradicting $\det P \equiv 1 \bmod \mf{m}^{\ltheta}$. 
	It follows that $\beta\gamma \in \I^{\ltheta}$ and  
	$\alpha\delta \equiv 1 \bmod \I^{\ltheta}$.
	Right multiplying $P$ by 
	the diagonal matrix of determinant $1$ with upper-left entry $\alpha^{-1}$,
	it is possible to choose both $\alpha$ and $\delta$ to be $1$ modulo $\I^{\ltheta}$.  
	Since $\alpha\beta \in B_{\theta}$ and $\gamma\delta \in C_{\theta}$,
	we must have $\beta \in B_{\theta}$ and $\gamma \in C_{\theta}$.		
	By Proposition \ref{rhoimage}, $P$ is then an element of $\rho_{\theta}(G_K)$.  It follows that 
	$\left( \begin{smallmatrix} -1& 0 \\ 0& 1 \end{smallmatrix} \right) \in 
	\rho_{\theta}(G_{\Q})$, as desired.
\end{proof}	

\begin{remark}
	One might ask if it is possible to use the same decomposition group $D_p$ 
	and element $\delta$ for all
	choices of $\theta$ in Theorem \ref{decomp}.  In fact, if $N = 1$, or if  
	$p$ does not divide $B_{1,\theta^{-1}}$ for any $\theta$ with $p \mid B_{1,\theta}$,
	then this follows quickly from an examination of the proofs of Proposition \ref{rhoimage}
	and Theorem \ref{decomp}.
	For $N=1$, the image of $\delta^{(p-1)/2}$ is the desired automorphism 
	for all $\theta$, which in this case are the odd powers of $\omega$.
	If $p$ does not divide any $B_{1,\theta^{-1}}$, then the analogue of
	Proposition \ref{rhoimage} holds for the full representation of $G_{\Q}$ on $\mc{Z}_N$,
	and the second method of proof in Theorem \ref{decomp} yields the result.
\end{remark}

\subsection{A comparison of Iwasawa and Hecke modules}

In \cite[Definition 4.1.17]{ohta-eich}, Ohta defines a
perfect $\zp[[1+p\zp]]$-bilinear pairing
\begin{equation*} 
	H^1_{\et}(N) \times H^1_{\et}(N) \to \zp[[1+p\zp]],
\end{equation*}
viewing $\zp[[1+p\zp]]$ as a subring of $\Lambda_N^{\mf{h}}$.  
We now give a slight modification of this.

Consider first the (canonical) twisted Poincar\'e duality pairing
$$
\label{Poincarer}
	(\,\cdot\,,\,\cdot\,)_r \colon 
	H^1(X_1^r(N),\zp)^{\ord} \times 
	H^1(X_1^r(N),\zp)^{\ord} \to \zp
$$
defined by the cup product
$$
	(x,y)_r = x\, \cup \, (w_{Np^r} (U_p^*)^r y),
$$ 
where $w_{Np^r}$ again denotes the Atkin-Lehner involution.  It is perfect and satisfies 
$$
	(T^*x,y)_r = (x,T^*y)_r
$$ 
for all $x, y \in H^1(X_1^r(N),\zp)^{\ord}$ and $T^* \in \mf{h}_r^*$.

\begin{proposition} \label{bigpairing}
	There exists a canonical perfect, $\Lambda_N^{\mf{h}}$-bilinear pairing		
	$$
		\langle\,\cdot\,,\,\cdot\, \rangle_N \colon 
		H^1(N) \times H^1(N) \to \Lambda_N^{\mf{h}}
	$$
	defined by the formula
	$$
		\langle x,y \rangle_N = \lim_{\substack{\leftarrow\\r}}
		\sum_{\substack{j=1\\(j,Np)=1}}^{Np^r-1} (x_r,
		\langle j^{-1} \rangle_r^*y_r)_r[j]_r
		\in \Lambda_N^{\mf{h}}
\label{PoincareN}
	$$	
	for $x = (x_r), y = (y_r) \in H^1(N)$ and satisfying 
	$$
		\langle T^*x,y \rangle_N = \langle x,T^*y \rangle_N
	$$ 
	for all $T^* \in \mf{h}^*$.
\end{proposition}
   
\begin{proof}
	Let $r \ge 1$.  
	The operator $w_{Np^{r+1}} U_p w_{Np^r}$ on $H_1(X_1^{r+1}(N);\zp)$
	is given by the sum
	$$
		\sum_{j=0}^{p-1} \left( \begin{matrix} 1 & 0 \\ jNp^r & 1\end{matrix} \right).
	$$
	Therefore, the natural (restriction) map from $H_1(X_1^r(N);\zp)$ is given by
	lifting and applying the operator
	$$
		w_{Np^{r+1}} U_p w_{Np^r} \sum_{k=0}^{p-1} \langle 1 + kNp^r \rangle.
	$$
	It follows then from Proposition \ref{poincare} and the compatibility of the comparison maps
	with restriction and Atkin-Lehner operators that the map
	$$
		\mr{Res} \colon H^1(X_1^r(N),\zp) \to
		H^1(X_1^{r+1}(N),\zp)
	$$
	that is identified with restriction on parabolic
	cohomology satisfies
	$$
		\mr{Res}(y_r) = w_{Np^{r+1}} U_p^*	w_{Np^r} 
		\sum_{k=0}^{p-1} \langle 1+kNp^r \rangle^*_{r+1} y_{r+1}.
	$$
	As the trace map commutes with $U_p^*$ and $w_{Np^r}$, it follows that
	$$
		\mr{Res}(w_{Np^r}(U_p^*)^r y_r) 
		= w_{Np^{r+1}} (U_p^*)^{r+1} \sum_{k=0}^{p-1}\langle 1+ kNp^r\rangle^*_{r+1} y_{r+1}
	$$
	and, therefore, that
	\begin{align*}
		(x_{r+1}, \sum_{k=0}^{p-1} \langle 1+kNp^r \rangle^*_{r+1} y_{r+1})_{r+1} = (x_r, y_r)_r, 
	\end{align*}
	Thus, the formula for $\langle x,y \rangle_N$ is well-defined.
	By definition, $\langle\,\cdot\,,\,\cdot\,\rangle_N$ 
	is $\Lambda_N^{\mf{h}}$-bilinear and satisfies the desired compatibility with
	the action of $\mf{h}^*$.  
	
	Since our pairing is $\Lambda_N^{\mf{h}}$-bilinear, its perfectness reduces
	to the question of the perfectness of the resulting pairings on
	eigenspaces
	$$
		\langle\,\cdot\,,\,\cdot\, \rangle_N^{\ltheta} \colon 
		H^1(N)^{\ltheta} \times 
		H^1(N)^{\ltheta} \to \Lambda_N^{\ltheta}
	$$
	for any character $\theta$ on $(\Z/Np\Z)^{\times}$.
	This is turn reduces to the perfectness of the pairing at level $Np$
	given by the projection of 
	$$
		\sum_{\substack{j=1\\(j,Np)=1}}^{Np-1} 
		(x_1,\langle j^{-1} \rangle_1^*y_1)_1 [j]_1
	$$
	to $\zp[(\Z/Np\Z)^{\times}]^{\ltheta}$
	for $x_1$, $y_1 \in (H^1(X_1^1(N),\zp)^{\ord})^{\ltheta}$.  This follows
	immediately from the perfectness of $(\,\cdot\,,\,\cdot\,)_1$.
\end{proof}

Using $\iota$, the pairing of Proposition \ref{bigpairing} allows us to define a $\Lambda_N^{\mf{h}}$-valued pairing on $H^1_{\et}(N)$, likewise denoted $\langle\,\cdot\,,\,\cdot\,\rangle_N$.  This is Galois equivariant with respect to the action of $\Gal(K/\Q)$ on $\Lambda_N^{\mf{h}}$ which, for the arithmetic Frobenius $\sigma_l$ attached to any prime $l \nmid Np$, is given by $(l[l])^{-1}$.
This follows from the fact that $w_{Np^r} \sigma_l \langle l \rangle^* = \sigma_l w_{Np^r}$ on $H^1_{\et}(X_1^r(N);\zp)$, together with Galois equivariance of the Poincar\'e duality pairing to $\zp(-1)$ (as in the proof of \cite[Corollary 4.2.8(ii)]{ohta-eich}).

\begin{proposition} \label{ohtapair}
	There exists a perfect, $\h/\I$-bilinear pairing
	$$ \label{finalpairing}
		\mc{Z}_N^-/\mc{Y}_N^- \times \mc{Y}_N^+/\I\mc{Y}_N^+ \to \h/\I,
	$$
	canonical up to the choice of $\iota$.
\end{proposition}

\begin{proof}
	We may consider the restriction of the pairing 
	$\langle\,\cdot\,,\,\cdot\,\rangle_N$ on \'etale cohomology 
	to a perfect pairing
	\begin{equation*}
		\langle\,\cdot\,,\,\cdot\,\rangle_N \colon
		\mc{Y}_N \times \mc{Y}_N \to \Lambda_N^{\mf{h}}
	\end{equation*}
	on Eisenstein parts.	
	Let $\theta$ denote an odd, primitive Dirichlet character of conductor $Np$.
	Restriction provides a perfect $\Lambda_N^{\ltheta}$-bilinear pairing
	\begin{equation*} \label{twistcup}
		\langle\,\cdot\,,\,\cdot\,\rangle_N^{\ltheta} \colon 
		\mc{Y}_N^{\ltheta} \times 
		\mc{Y}_N^{\ltheta}
		\to \Lambda_N^{\ltheta},
	\end{equation*}
	satisfying the same Hecke compatibility as $\langle\,\cdot\,,\,\cdot\, \rangle_N$.
		
	Since $\mc{Z}_N^{\ltheta}$ and $\mc{Y}_N^{\langle \theta
	\rangle}$ are both free of the same 
	$\Lambda_N^{\ltheta}$-rank with $\mc{Y}_N^{\langle \theta
	\rangle} \subset \mc{Z}_N^{\ltheta}$,
	we may extend $\langle\,\cdot\,,\,\cdot\,\rangle_N^{\ltheta}$
	uniquely to a pairing
	\begin{equation} \label{extend}
		\mc{Z}_N^{\ltheta} 
		\times \mc{Y}_N^{\ltheta} \to \mf{L}_N^{\ltheta}
	\end{equation}
	to the quotient field of $\Lambda_N^{\ltheta}$.
	The aforementioned Galois equivariance implies that 
	$(\mc{Z}_N^{\ltheta})^{\pm}$ pairs trivially with 	$(\mc{Y}_N^{\ltheta})^{\pm}$.
	Note that $(\mc{Y}_N^{\ltheta})^+/\I (\mc{Y}_N^{\ltheta})^+$ and 
	$(\mc{Z}_N^{\ltheta})^-/(\mc{Y}_N^{\ltheta})^-$ are both isomorphic to 
	$(\mf{h}^*/\I)^{\ltheta}$ as Hecke modules. 
	Reducing modulo $\Lambda_N^{\ltheta}$ and taking the direct sum over
	a set of representatives for the odd classes in $\Sigma_{Np}$, we 
	finally obtain our pairing.
\end{proof}

\begin{lemma} \label{xihomology}
	Let $r \ge 1$, and let $u$ and $v$ be positive integers not 
	divisible by $Np^r$ that satisfy $(u,v,Np) = 1$.
	Then $e_r[u:v]_r \in H_1(X_1^r(N);\zp)_{\mf{m}}$.
\end{lemma}

\begin{proof}
	As in Lemma \ref{inhomology}, this reduces to showing 
	that any cusp $\binom{a}{bM}_r$ with $M$ a nontrivial divisor of $Np^r$ 
	has trivial image in $H_0(C_1^r(N);\zp)_{\mf{M}}$.  
	We follow the argument of \cite[Proposition 4.3.4]{ohta-ord}.
	Suppose $l$ is a prime dividing $M$,
	and let $s$ be such that $l^s$ exactly divides $Np^r/M$. 
	Let $t > s$ be such that $l^{t-s} \equiv 1 \bmod P$, where $P$
	denotes the prime-to-$l$ part of $Np^r$.
	Then
	$$
		U_l^t\binom{a}{bM}_r = 
		\sum_{i=0}^{l^t-1} \binom{a+bMi}{l^tbM}_r
		= l^{t-s} \sum_{i=0}^{l^s-1} \binom{a+bMi}{l^sbM}_r
		= l^{t-s}U_l^s\binom{a}{bM}_r.
	$$
	But $U_l$ acts as $1$ on the free $\zp$-module $H_0(C_1^r(N);\zp)_{\mf{M}}$ 
	(see \cite[Theorem 2.3.6]{ohta-cong}), so we must
	have $l^{t-s}=1$ in $\zp$ or $\binom{a}{bM}_r = 0$.  Clearly,
	the former is impossible.
\end{proof}

We have the following immediate corollary.

\begin{corollary}
	Let $u \in \Z[\frac{1}{p}]$ be nonzero, let $v \in \Z$ be prime to $p$,
	and suppose that $(u,v,N)\Z[\frac{1}{p}] = \Z[\frac{1}{p}]$.  Then
	$\xi(u:v) \in H_1(N)_{\mf{m}}$.
\end{corollary}

\begin{lemma} \label{xi0}
	We have that $\I\mc{Z}_N^- \subseteq \mc{Y}_N^-$, and
	the image of $\xi(0:1) \in \mc{H}_1(N)^+$ generates
	$\mc{Z}_N^-/\mc{Y}_N^-$ as an $\h/\I$-module.
\end{lemma}

\begin{proof}
	As in the proof of Lemma \ref{usefulfact}, we have
	$$
		\mc{Z}_N/\mc{Y}_N \cong \mc{Z}_N^-/\mc{Y}_N^- \cong \h/\I,
	$$ 
 	and $\h/\I$ is canonically isomorphic to a quotient of $\Lambda_N^{\mf{h}}$.
	The first statement follows.
	Note that $\mc{H}_1(N)^{\pm}$ is isomorphic to $\mc{H}^1(N)^{\mp}$, 
	since complex conjugation acts on $H^2_c(X_1^r(N);\zp)$ as $-1$.  Hence, the image of 
	$\xi(0:1) \in \mc{H}_1(N)^+$ in $\mc{H}^1(N)$ lies in $\mc{H}^1(N)^-$.  That its image in 
	$\mc{Z}_N^-/\mc{Y}_N^-$ is a generator follows the definition of the
	congruence module and the proof of 
	\cite[Theorem 2.3.6]{ohta-cong}, since it is shown there that the projection
	of the cusp $\binom{0}{1}_r$ to the Eisenstein component of 
	$\tilde{H}_0(C_1^r(N);\zp)$ generates it as a Hecke module (and 
	we know that the image of $\infty$ is trivial).
\end{proof}

The pairing of Proposition \ref{ohtapair}
induces an isomorphism of $\h/\I$-modules:
\begin{equation} \label{pairingiso}
	\mc{Z}_N^-/\mc{Y}_N^- \xrightarrow{\sim} 
	\Hom_{\h}(\mc{Y}_N^+/\I\mc{Y}_N^+,\h/\I).
\end{equation}
We therefore have the following corollary.

\begin{corollary}
\label{Pi}
	The map $\Pi$ induced by applying \eqref{pairingiso} to the image of
	$\xi(0:1) \in \mc{H}_1(N)^+$ in $\mc{Z}_N^-/\mc{Y}_N^-$
	generates
	$
		\Hom_{\h}(\mc{Y}_N^+/\I\mc{Y}_N^+,\h/\I)
	$
	as a Hecke module and is canonical up to the choice of $\iota$.
\end{corollary}

Let 
\label{XK}
$X_K$ denote the Galois group of the 
maximal unramified abelian pro-$p$ extension of $K$.
Now, we compare the $\Lambda_N$-modules of interest.

\begin{proposition} \label{1stisom}
	We have a homomorphism
\label{phi1'}
	$$
		\phi_1' \colon X_K^{\circ} \to
		\mc{Y}_N^-/\I \mc{Y}_N^-
	$$
	of $\Lambda_N$-modules $($under Galois$)$, 
	canonical up to the choice of $\iota$,
	which is an isomorphism in its $(\theta^{-1})$-eigenspace for
	$\theta$ odd and primitive if $p \nmid B_{1,\theta^{-1}}$.
\end{proposition}

\begin{proof}
	The Galois action on $\mc{Y}_N$ provides a map
	$$
		b \colon G_{\Q} \to \Hom_{\h}(\mc{Y}_N^+,\mc{Y}_N^-)
	$$
	that, by Theorem \ref{decomp} and \cite[Theorem 3.3.12]{ohta-ord2} 
	(see also \cite[Proposition 1.8.2]{mw1}),
	induces a homomorphism
	$$
		\bar{b} \colon X_K^{\circ} \to \Hom_{\h}(\mc{Y}_N^+/\I \mc{Y}_N^+,
		\mc{Y}_N^-/\I \mc{Y}_N^-)
	$$
	of Galois $\Lambda_N$-modules.
	Since
	\begin{equation*} \label{homgptens}
		\Hom_{\h}(\mc{Y}_N^+/\I\mc{Y}_N^+,\mc{Y}_N^-/\I\mc{Y}_N^-) \cong
		\Hom_{\h}(\mc{Y}_N^+/\I\mc{Y}_N^+,\h/\I) \otimes_{\h/\I} 
		\mc{Y}_N^-/\I \mc{Y}_N^-,
	\end{equation*}
	we may define $\phi_1'$ by $\bar{b}(\sigma) = \Pi \otimes \phi_1'(\sigma)$
	for $\sigma \in X_K^{\circ}$ and the generator $\Pi$ defined above.
	
	Let $B_N$ denote the Hecke submodule of $\mc{Y}_N^-$ generated by
	the images of elements in the image of $b$.
	As $\mc{Y}_N^+$ is isomorphic to $\h$ as an $\h$-module, $B_N^{\ltheta}$
	is isomorphic to the $(\theta\omega^{-1})$-eigenspace of the image of $b$ for the action
	of the adjoint diamond operators in $(\Z/Np\Z)^{\times}$.  That $\phi_1'$ 	
	is an isomorphism in its $(\theta^{-1})$-eigenspace under Galois then follows from
	\cite[(5.3.18) and (5.3.20)]{ohta-ord} (and \cite[Section 3.2]{ohta-ord2}) whenever
	$B_N^{\ltheta} = (\mc{Y}_N^-)^{\ltheta}$. 
	If $p \nmid B_{1,\theta^{-1}}$, then $(\mc{Z}_N^-)^{\ltheta}$ 
	is free of rank $1$ over $(\h)^{\ltheta}$ 
	by \cite[(3.4.7)]{ohta-comp}.  Thus, the fact
	that $\mc{Z}_N^-/\mc{Y}_N^- \cong \h/\I$ implies that 
	$(\mc{Y}_N^-)^{\ltheta} = (\I\mc{Z}_N^-)^{\ltheta}$. 
	Combining this with \cite[(3.4.10)]{ohta-comp}, which tells us that
	$B_N^{\ltheta} = (\I\mc{Z}_N^-)^{\ltheta}$, 
	we obtain the final part of the proposition.
\end{proof}

For the purpose of formulating our conjectures,
we will use an ill-defined modification of $\phi_1'$ throughout. That is, we set 
\label{phi1}
$$
	\phi_1 = c_N\phi_1'
$$ 
for a fixed unit 
\label{cN}
$c_N \in \Lambda_N^{\circ}$, independent of $\iota$, that makes Conjecture \ref{symbolconj} below true.  
Though we do not write $c_N$ directly into the statements of our conjectures, one should, of course, still understand its existence to be a part of them.

\subsection{Inverse limits of cup products and modular symbols}

First, we show that for the purposes of considering the primitive part
of second cohomology group, it suffices to restrict to the primitive part of 
the maximal unramified abelian pro-$p$ extension.

\begin{lemma} \label{h2isom}
	The canonical homomorphism $X_K^{\circ} \to H_S^2(K,\zp(1))^{\circ}$ is an 
	isomorphism.
\end{lemma}

\begin{proof}
	Note first that, as seen in the proof of \cite[Lemma 3.4]{me-paireis}, 
	we have $X_K^{\circ} = X_{K,S}^{\circ}$.
	If $l$ is a prime dividing $Np$, then 
	the part of the direct sum in \eqref{h2seq} arising from primes over $l$
	has a trivial $(\Z/l\Z)^{\times}$-action since
	there is a unique prime above $l$ in 
	$\Q(\mu_l)$.  But this means, in particular,
	that the primitive part of the direct sum in \eqref{h2seq} must be
	trivial.
\end{proof}

We will let 
\label{cupKScirc}
$(\,\cdot\,,\,\cdot\,)_{K,S}^{\circ}$ denote the projection of the pairing $(\,\cdot\,,\,\cdot\,)_{K,S}$ to $X_K^{\circ}(1)$.  For $v$ prime to $p$, we let 
\label{1zeta}
$1-\zeta^v$ denote the norm compatible sequence of elements $1-\zeta_{Np^r}^v \in \Q(\mu_{Np^r})$.
We use $\Upsilon_K$ to denote the map
$$
\label{UpsilonK}
	\Upsilon_K \colon X_K^{\circ}(1) \to 
	(\mc{Y}_N^-/\I\mc{Y}_N^-)(1)
$$
that is the Tate twist of $\phi_1$.

Recall that $\mc{H}_1(N)_{\mf{m}}^+ \cong \mc{Y}_N^-(1)$ canonically
up to the choice of $\iota$.  
For $u \in \Z[\frac{1}{p}]$ and $v \in \Z$ prime to $p$ with
$(u,v,N)\Z[\frac{1}{p}] = \Z[\frac{1}{p}]$, we let
\label{barxi}
$\bar{\xi}(u:v)$ denote the image of $\xi(u:v)^+$ in 
$(\mc{Y}_N^-/\I\mc{Y}_N^-)(1)$.

We may now phrase the first form of our conjecture as follows. 

\begin{conjecture} \label{symbolconj}
	For any $s \ge 0$ and all nonzero $u$ and $v \in \Z$ prime to $p$ with
	$(u,v,N) = 1$ and $u$ not divisible by $Np^s$,
	we have 
	$$
		\Upsilon_K((1-\zeta_{Np^s}^u,1-\zeta^v)_{K,S}^{\circ}) =
		\bar{\xi}(p^{-s}u:v).
	$$ 
\end{conjecture}

We verify the independence of Conjecture \ref{symbolconj} from the choice of the complex embedding $\iota$.

\begin{proposition}
	The validity of Conjecture \ref{symbolconj} is independent of the
	choice of $\iota$.
\end{proposition}

\begin{proof}
	Let us choose a second complex embedding $\iota'$ 
	of the form $\iota' = \iota \circ \sigma^{-1}$ for some $\sigma \in G_{\Q}$.
	Let $\mf{Y}_N^{\pm}$
	denote the $(\pm 1)$-eigenspaces of $\mc{Y}_N$ 
	under the complex conjugation determined by $\iota'$.  
	Let $\Upsilon_K'$, $\phi_1'$, and $\bar{b}'$ denote
	the maps arising from Proposition \ref{1stisom}
	with the embedding $\iota'$.  Let $\Pi'$ denote the map
	defined as $\Pi$ using $\iota'$.  We use $\bar{\xi}'$ in denoting the
	symbols defined using $\iota'$ and corresponding to those denoted with 
	$\bar{\xi}$.
	
	Let us first consider the map $\phi_1$.  
	Note that $\mf{Y}_N^{\pm} = \sigma(\mc{Y}_N^{\pm})$ and, by construction, we have
	$$
		\Pi'(y') = \Pi(\sigma^{-1}y')
	$$
	for $y' \in \mf{Y}_N^+/\I \mf{Y}_N^+$.
	We therefore have a commutative diagram
	$$
		\SelectTips{cm}{} 
		\xymatrix{
		X_K^{(\theta^{-1})} \ar[r]^-{\bar{b}} 
		\ar[d]_{\tau \mapsto \sigma \tau\sigma^{-1}} &
		\Hom_{\h}(\mc{Y}_N^+/\I\mc{Y}_N^+,\mc{Y}_N^-/\I \mc{Y}_N^-)
		\ar[d]^{f \mapsto \sigma\circ f \circ\sigma^{-1}} 
		\ar[r]^-{\Pi \otimes 1} & \mc{Y}_N^-/\I \mc{Y}_N^- 
		\ar[d]^{y \mapsto \sigma y}
		 \\
		X_K^{(\theta^{-1})} \ar[r]^-{\bar{b'}} & 
		\Hom_{\h}(\mf{Y}_N^+/\I \mf{Y}_N^+,
		\mf{Y}_N^-/\I \mf{Y}_N^-) \ar[r]^-{\Pi' \otimes 1} &
		\mf{Y}_N^-/\I \mf{Y}_N^-.}
	$$
	In other words, we have
	\begin{equation} \label{map1}
		\phi_1'(\sigma\tau\sigma^{-1}) = \sigma\phi_1(\tau).
	\end{equation}
	for $\tau \in X_K^{(\theta^{-1})}$.
		
	Next, note that the change of embedding takes $1-\zeta^v$ to 
	$\sigma(1-\zeta^v)$ and $1-\zeta_{Np^s}^u$ to $\sigma(1-\zeta_{Np^s}^u)$.  
	Using \eqref{map1} (applied to $\Upsilon_K$) and the
	Galois equivariance of $(\,\cdot\,,\,\cdot\,)_{K,S}^{\circ}$, we see that
	\begin{equation} \label{imageinvar}
		\Upsilon'_K((\sigma(1-\zeta_{Np^s}^u),\sigma(1-\zeta^v))_{K,S}^{\circ})
		= \sigma\Upsilon_K((1-\zeta_{Np^s}^u,1-\zeta^v)_{K,S}^{\circ}).
	\end{equation}
	
	On the other hand, we have a map $\alpha$ that is the isomorphism
	$$
		\mc{H}_1(N) \xrightarrow{\Phi} \mc{H}_1^{\et}(N) \xrightarrow{\sim} 
		\mc{H}^1_{\et}(N)(1)
	$$
	and the analogous map $\alpha'$ defined using $\iota'$.
	One sees immediately that $\alpha' = \sigma \circ \alpha$.
	It follows that
	\begin{equation} \label{Lfninvar}
		\bar{\xi}'\left(p^{-s}u:v\right) =
		\sigma\bar{\xi}\left(p^{-s}u:v\right).
	\end{equation}
	Comparing \eqref{Lfninvar}
	with \eqref{imageinvar},
	we see that if Conjecture \ref{mainconj}
	holds with $\iota$, it must also hold with $\iota'$. 
\end{proof}

\section{The view from finite level}

\subsection{Cup products and modular symbols} \label{cpms}
We now consider the implications of Conjecture \ref{symbolconj} at finite
level.  For now, we focus on weight $2$.  Let $r \ge 1$, and set 
\label{Fr}
$F_r = F(\mu_{p^r})$.  Let 
\label{Yr}
$Y_r = H^1_{\et}(X_1^r(N),\zp)_{\mf{m}}$, and let 
\label{Ir}
$I_r$ denote
the image of $\mc{I}$ in $\mf{h}_r^*$.
Let 
$$
\label{H2circ}
	H^2_{\cts}(G_{F_r,S},\zp(2))^{\circ} = \bigoplus_{\substack{(\chi) \in \Sigma_{Np}
	\\ \chi \text{ odd}}} H^2_{\cts}(G_{F_r,S},\zp(2))^{(\chi\omega)}.
$$
Note that this differs slightly from our previous version of $A^{\circ}$
due to the twist in the cohomology group.

\begin{lemma} \label{weight2map}
	For each $r \ge 1$, there exists a map 
	$$
\label{nur}
		\nu_r \colon H^2_{\cts}(G_{F_r,S},\zp(2))^{\circ} 
		\to (Y_r^-/I_rY_r^-)(1),
	$$
	canonical up to the choice of $\iota$,
	that is an isomorphism in its $(\omega\theta^{-1})$-eigenspace if
	$p \nmid B_{1,\theta^{-1}}$.
\end{lemma}

\begin{proof}
	We construct $\nu_r$ out of $\Upsilon_K$.
	By Lemma \ref{h2isom},
	we have 
	$$
		H^2_S(K,\zp(2))^{\circ} \cong X_K^{\circ}(1).
	$$ 
	Since $G_{F,S}$ has $p$-cohomological
	dimension $2$, 
	corestriction then defines an isomorphism
	\begin{equation} \label{YKcoinv}
		X_K^{\circ}(1)_{\Gamma_r} \xrightarrow{\sim} 
		H^2_{\cts}(G_{F_r,S},\zp(2))^{\circ}
	\end{equation}
	with $\Gamma_r = \Gal(K/F_r)$.
	
	Set $\omega_r = (\langle 1+p \rangle^*)^{p^{r-1}} 
	- 1 \in \mf{h}^*$.  By \cite[Theorem 1.2]{hida}, we have that
	$\mf{h}^*/\omega_r\mf{h}^* \cong \mf{h}^*_r$,  
	and by \cite[Theorem 1.4.3]{ohta-eich}, we have
	$\mc{Y}_N/\omega_r\mc{Y}_N \cong Y_r$.
	The Galois element $\sigma_j$ corresponding to $j \in 
	\Z_{p,N}^{\times}$ acts on $\mc{Y}_N^-/\I\mc{Y}_N^-$ as $(\chi(j)\langle j \rangle^*)^{-1}$,
	where $\chi$ is the $p$-adic cyclotomic
	character.  Thus, corestriction provides an isomorphism
	$$
		(\mc{Y}_N^-/\I\mc{Y}_N^-)(1)_{\Gamma_r} \cong (Y_r^-/I_rY_r^-)(1).  
	$$ 
	We take $\nu_r$ to be the map arising
	from $\Upsilon_K$ on $\Gamma_r$-coinvariants.  The final statement
	now follows from the final statement of Proposition \ref{1stisom}.
\end{proof} 

Let us denote the pairing induced from $(\,\cdot\,,\,\cdot\,)_{F_r,S}$ via 
projection to $H^2_{\cts}(G_{F_r,S},\zp(2))^{\circ}$
by 
\label{cupFr}
$(\,\cdot\,,\,\cdot\,)_{F_r,S}^{\circ}$.
For $u, v \in \Z$ not divisible by $Np^r$ and with $(u,v,Np^r) = 1$, we let
$\bar{\xi}_r(u:v)$ denote the image of $\xi_r(u:v)^+$ in $Y_r^-/I_rY_r^-$
(which depends only upon $u$ and $v$ modulo $Np^r$).
We now state an analogue of Conjecture \ref{symbolconj} at the finite level.

\begin{conjecture} \label{finsymbconj}
	Let $r \ge 1$.  
	Suppose that $u$ and $v$ are positive integers not divisible by $Np^r$
	with $(u,v,Np) = 1$.  
	Then we have
	$$
		\nu_r((1-\zeta_{Np^r}^u,1-\zeta_{Np^r}^v)_{F_r,S}^{\circ})
		= \bar{\xi}_r(u:v).
	$$
\end{conjecture}

In fact, this conjecture is equivalent to Conjecture \ref{symbolconj}.

\begin{proposition} \label{equivconjs}
	Conjecture \ref{symbolconj} and Conjecture \ref{finsymbconj}
	are equivalent.
\end{proposition}

\begin{proof}
	Let $u$, $v$, and $s$ be as in Conjecture \ref{symbolconj}.
	The corestriction map yielding
	\eqref{YKcoinv} takes $(1-\zeta_{Np^s}^u,1-\zeta^v)_{K,S}^{\circ}$
	to $(1-\zeta_{Np^s}^u,1-\zeta_{Np^r}^v)_{F_r,S}^{\circ}$, 
	and the map 
	$$\mc{H}_1(N) \to H_1(X_1(Np^r);\zp)^{\ord}$$ 
	takes $\xi(p^{-s}u:v)$ to $\xi_r(p^{r-s}u:v)$ for $r \ge s$.  
	Since in each of the two cases,
	the former object is the inverse limit of the latter objects,
	we have both implications.
\end{proof}

Suppose that $t$ is a positive divisor of $Np^r$ for some $r$ and $u$ and $v$ are positive nonmultiples of $Np^r$ with
$(tu,v,Np) = 1$, and set $Q = Np^r/t$.  We also assume that $u$ is not a multiple of $Q$.  Since $U_t-1 \in \I$, the equation \eqref{maninsum} yields
immediately that
\begin{equation} \label{maninredsum}
	\sum_{k=0}^{t-1}\bar{\xi}_r(u+kQ:v) = \bar{\xi}_r(tu:v).
\end{equation}
On the other hand,
$$
	\sum_{k=0}^{t-1}(1-\zeta_{Np^r}^{u+kQ},1-\zeta_{Np^r}^v)_{F_r,S}
	= (1-\zeta_Q^u,1-\zeta_{Np^r}^v)_{F_r,S}.
$$
In particular, Conjecture \ref{finsymbconj} is compatible with these relations.

\subsection{Image of the cup product pairing}

We have the following generalization of a conjecture of McCallum and the author's \cite[Conjecture 5.3]{mcs}, originally given in the case $N=1$.

\begin{conjecture} \label{pairconj}
	The span of the image of $(\,\cdot\,,\,\cdot\,)_{F_r,S}^{\circ}$ is
	$H^2_{\cts}(G_{F_r,S},\zp(2))^{\circ}$.
\end{conjecture}

We require the following lemma.

\begin{lemma} \label{allmaningen}
	The images of the 
	symbols $[u:v]^+_r$ for nonzero $u, v \in \Z/Np^r\Z$ with $(u,v) = (1)$
	together generate $H_1(X_1^r(N);\zp)^+_{\mf{m}}$ as a $\zp$-module.
\end{lemma}

\begin{proof}
	Lemma \ref{xihomology} implies that such a 
	$[u:v]_r$ lies in $H_1(X_1^r(N);\zp)_{\mf{m}}$ since $u, v \neq 0$.
	Furthermore, \cite[Theorem 2.3.6]{ohta-cong} implies
	that $\tilde{H}_0(C_1^r(N);\zp)^+_{\mf{M}}$
	is freely generated as a module over the image of $\zp[(\Z/Np^r\Z)^{\times}]$ in $(\mf{H}_r)_{\mf{M}}$
	by the image of $[0:1]_r$.  
	Since the $\zp[(\Z/Np^r\Z)^{\times}]$-span 
	of $[0:1]_r$ in $H_1(X_1^r(N),C_1^r(N);\zp)_{\mf{M}}$
	contains the $[0:w]_r$ with $1 \le w < Np^r$ and $(w,Np) = 1$,
	the exact sequence \eqref{homologyseq} yields the result.
\end{proof}

We now see that Conjecture \ref{finsymbconj} implies much of Conjecture
\ref{pairconj}.

\begin{proposition}
	Conjecture \ref{finsymbconj} implies that the span of 
	the image of $(\,\cdot\,,\,\cdot\,)_{F_r,S}^{\circ}$ contains
	$H^2_{\cts}(G_{F_r,S},\zp(2))^{(\omega\theta^{-1})}$
	for all primitive odd $\theta$ with $p \nmid B_{1,\theta^{-1}}$.
\end{proposition}

\begin{proof} 
	Since $p \nmid B_{1,\theta^{-1}}$,
	Proposition \ref{1stisom} implies that $\nu_r$ is an isomorphism.
	Lemma \ref{allmaningen} and Conjecture \ref{finsymbconj} then imply 
	that the images of the $(1-\zeta_{Np^r}^u,1-\zeta_{Np^r}^v)_{F_r,S}$
	generate the $(\omega\theta^{-1})$-eigenspace of 
	$H^2_{\cts}(G_{F_r,S},\zp(2))$, as desired.
\end{proof}

\subsection{A map in the other direction} \label{mapback}

The comparison between the two sides of Conjecture \ref{finsymbconj} is 
perhaps seen more naturally in the opposite direction.  We begin by examining relations among values of the cup product pairings on cyclotomic $S$-units.
We will find relations analogous to the relations \eqref{maninrels1}-\eqref{maninrels5} on Manin symbols.

Recall that $(x,1-x)_{F_r,S} = 0$ if $x$ and $1-x$ are both $S$-units in $F_r$
\cite[Corollary 2.6]{mcs}.  
Note that $(\zeta_{Np^r},\zeta_{Np^r})_{F_r,S} = 0$ by antisymmetry
of the cup product. Since
$\mc{E}_{F_r} \cong \mc{E}_{F_r}^+ \oplus \mu_{p^r}$ and $H^2_{\cts}(G_{F_r,S},\zp(2))^{\circ}$ has a trivial action of $-1$, Galois equivariance of the
cup product pairing implies that 
$(\zeta_{Np^r},x)_{F_r,S}^{\circ} = 0$ for all $x \in \mc{E}_{F_r}$.

Suppose that $u$ and $v$ are integers that are not divisible by $Np^r$.
Since 
$$
	1-\zeta_{Np^r}^u = -\zeta_{Np^r}^u(1-\zeta_{Np^r}^{-u}),
$$ 
we have
\begin{equation} \label{minus}
	(1-\zeta_{Np^r}^u,1-\zeta_{Np^r}^v)_{F_r,S}^{\circ}
	= (1-\zeta_{Np^r}^{u},1-\zeta_{Np^r}^{-v})_{F_r,S}^{\circ}
	= (1-\zeta_{Np^r}^{-u},1-\zeta_{Np^r}^{-v})_{F_r,S}^{\circ}.
\end{equation}
Antisymmetry of the cup product yields
\begin{equation} \label{antisymm}
	(1-\zeta_{Np^r}^u,1-\zeta_{Np^r}^v)_{F_r,S} +
	(1-\zeta_{Np^r}^v,1-\zeta_{Np^r}^u)_{F_r,S} = 0.
\end{equation}
Additionally, if $u+v$ is not divisible by $Np^r$, then the identity
$$
	1-\zeta_{Np^r}^u + 
	\zeta_{Np^r}^u(1-\zeta_{Np^r}^v) = 1-\zeta_{Np^r}^{u+v},
$$
implies
\begin{equation} \label{triple}
	(1-\zeta_{Np^r}^u,1-\zeta_{Np^r}^v)_{F_r,S}^{\circ} = 
	(1-\zeta_{Np^r}^u,1-\zeta_{Np^r}^{u+v})_{F_r,S}^{\circ} +
	(1-\zeta_{Np^r}^{u+v},1-\zeta_{Np^r}^v)_{F_r,S}^{\circ}.
\end{equation}
Finally, Galois equivariance tells us that, for any $j \in \Z$ prime
to $Np$, we have
\begin{equation} \label{sigma}
	\sigma_j(1-\zeta_{Np^r}^u,1-\zeta_{Np^r}^v)_{F_r,S} 
	= (1-\zeta_{Np^r}^{ju},1-\zeta_{Np^r}^{jv})_{F_r,S},
\end{equation}
where $\sigma_j \in \Gal(K/\Q)$ satisifes $\sigma_j(\zeta_{Np^r})
= \zeta_{Np^r}^j$.

\begin{proposition} \label{otherdir}
	There exists a homomorphism
	$$
\label{varpir}
		\varpi_r \colon H_1(X_1^r(N),C_1^r(N);\zp)^+ 
		\to H^2_{\cts}(G_{F_r,S},\zp(2))^{\circ}
	$$
	satisfying  $\varpi_r \circ \langle j \rangle_r = \sigma_j^{-1} \circ \varpi_r$
	for all $j$ prime to $Np$ and such that
	$\varpi_r([1:0]_r^+) = 0$ and
	$$
		\varpi_r([u:v]_r^+) = 
		(1-\zeta_{Np^r}^u,1-\zeta_{Np^r}^v)_{F_r,S}^{\circ}
	$$
	for $u, v \in \Z$ not divisible by $Np^r$ with $(u,v,Np) = 1$.
\end{proposition}

\begin{proof}
	Compare relations \eqref{maninrels1}, \eqref{maninrels3}, and
	\eqref{maninrels5} with \eqref{minus} and \eqref{antisymm}, relation 
	\eqref{maninrels2} with \eqref{triple}, and relation \eqref{maninrels4}
	with \eqref{sigma}.  Since the
	Manin symbols generate the module
	$H_1(X_1^r(N),C_1^r(N);\zp)^+$ and the relations
	\eqref{maninrels1}-\eqref{maninrels5} give a presentation of it,	
	we need only remark that $\varpi_r$ behaves well with respect to these
	relations in the case $v = 0$.  This is obvious: for instance, for
	relation \eqref{maninrels2}, we have
	$$
		\varpi_r([u:0]_r^+) + \varpi_r([u:u]_r^+) + \varpi_r([0:u]^+) = 0,
	$$
	as $[u:u]_r^+ = 0$.
\end{proof}

We fully expect that the restriction of $\varpi_r$ to $H_1(X_1^r(N);\zp)^+$ is
Eisenstein.  

\begin{conjecture} \label{oppconj}
	The restriction of  $\varpi_r$ to $H_1(X_1^r(N);\zp)^+$ 
	satisfies 
	$$
		\varpi_r(T_lx) = (1+l\varepsilon(\langle l \rangle))
		\varpi_r(x)
	$$
	for $l \nmid Np$ and
	$$
		\varpi_r(U_lx) = \varpi_r(x)
	$$
	for $l \mid Np$, for all $x \in H_1(X_1^r(N);\zp)^+$.
\end{conjecture}

One can check directly using relations of McCallum and the author in Milnor $K_2$ of $\mc{O}_{F_r,S}$ (e.g., \cite[Section 5]{mcs} for $T_2$) that it is Eisenstein with respect to the operators $T_2$ if $2 \nmid N$ and $T_3$ if $3 \nmid N$.
A slight variant of the map $\varpi_r$ and this fact have been discovered independently by C. Busuioc, and we refer the reader to \cite{busuioc} for the latter (in the case $N=1$), as our proof is very similar.  The analogue of Conjecture \ref{oppconj} is also discussed in \cite{busuioc}.  (We remark that in \cite[Section 9]{busuioc}, Vandiver's conjecture at $p$ should be assumed for the conjecture to hold in the form stated.)

\begin{remark}
	When taken together with Conjectures \ref{finsymbconj} and \ref{pairconj}
	(noting Lemma \ref{allmaningen}), Proposition \ref{otherdir}
	forces $\nu_r$ and the map that $\varpi_r$ induces on $(Y_r^-/\mc{I}Y_r^-)(1)$
	to be inverse isomorphisms.  It follows that $\Upsilon_K$ would also
	be an isomorphism, or equivalently, that $B_N$ in the proof of Proposition
	\ref{1stisom} would equal $\mc{Y}_N^-$.  While rather natural, this is also 
	a remarkably strong statement, which does give us some pause.
\end{remark}

\section{Main form of the conjecture}

\subsection{Modified two-variable $p$-adic $L$-functions} \label{modified}

Let $M$ be a positive divisor of $N$.  We now consider modifications of our two-variable
$p$-adic $L$-functions
$\mc{L}_{N,M}$.  We view $\zp[[\Z_{p,N}]]$ as a continuous module over $\Lambda_N$ via left multiplication.  We again denote the element of $\zp[[\Z_{p,N}]]$ corresponding to $j \in \Z_{p,N}$ by 
\label{[j]2}
$[j]$.
Define 
\label{LambdaNstar}
$\Lambda_N^{\star}$ to be the quotient of $\zp[[\Z_{p,N}]]$ by the
$\Lambda_N$-submodule generated by $[0]$.
Define 
$\Z_{p,N}^{\star}$ to be the set of nonzero elements in $\Z_{p,N}$.
We also write 
\label{ZpNstar}
$\Lambda_N^{\star} = \zp[[\Z_{p,N}^{\star}]]$. 

We now construct modified versions of our $L$-functions.
For any $M \ge 1$ dividing $N$, let us set
\begin{equation}
\label{LNMstar}
	\mc{L}_{N,M}^{\star} = 
	\lim_{\leftarrow} \sum_{\substack{j=1\\(j,M) = 1}}^{Np^r-1} 
	U_p^{-r}\xi_r(j:M)
	\otimes [j]_r \in \mc{H}_1(N) \cotimes \Lambda_N^{\star}.
\end{equation}  
It is worth noting that this is well-defined.  The proof is similar to the case
of $\mc{L}_N$, with one additional detail.  
That is, we view 
\label{[j]r2}
$[j]_r$ in the
$r$th term of the inverse limit in $\eqref{LNMstar}$ as an element of
$\zp[(\Z/Np^r\Z)^{\star}]$, which we define to be the quotient of 
$\zp[\Z/Np^r\Z]$
by the $\zp[(\Z/Np^r\Z)^{\times}]$-submodule generated by $[0]_r$.
The natural map
$$
	\zp[(\Z/Np^s\Z)^{\star}] \to \zp[(\Z/Np^r\Z)^{\star}]
$$
for $s \ge r$ now takes $[j]_s$ to $[j]_r$, the latter of which is $0$ if $j \equiv 0 \bmod Np^r$.  The rest is then the same as before.  Note that we use $(\Z/Np^r\Z)^{\star}$
to denote the set of nonzero elements in $\Z/Np^r\Z$.

Although $\mc{L}_{N,M}^{\star}$ will not always lie in 
$H_1(N) \cotimes \Lambda_N^{\star}$, its localization
in the Eisenstein part of homology does.  The following is an immediate
corollary of Lemma \ref{xihomology}.

\begin{corollary} \label{LN-}
	The modified $L$-function $\mc{L}_{N,M}^{\star}$ lies in 
	$H_1(N)_{\mf{m}} \cotimes \Lambda_N^{\star}$.
\end{corollary}

Finally, we remark that $\mc{L}_{N,M}^{\star}$ also specializes
to integrals with respect to $\lambda_{N,M}$.  We extend 
\label{lambdaNM2}
$\lambda_{N,M}$ to a measure on $\Z_{p,N}$ by setting
$$
	\lambda_{N,M}(a+Np^r\Z_{p,N}) = 
	\begin{cases}
	0 & \text{if }(a,M) \neq 1 \\
	U_p^{-r}\xi(a:M) & \text{otherwise}.
	\end{cases}
$$  
Let $\chi \colon \Z_{p,N} \to \overline{\qp}$ be a congruence function (i.e., a uniform
limit of congruence functions of finite period, necessarily satisfying
$\chi(0) = 0$).
Consider the induced map 
$$
\label{tildechi2}
	\tilde{\chi} \colon \mc{H}_1(N) \cotimes \Lambda_N^{\star} \to
	\mc{H}_1(N) \otimes_{\zp} \qpbar.
$$
We then have
\begin{equation}
	\tilde{\chi}(\mc{L}_{N,M}^{\star}) = 
	\int_{\Z_{p,N}} 
	\chi\lambda_{N,M} \in \mc{H}_1(N) \otimes_{\zp} \qpbar.
\end{equation}

\subsection{The $\zp$-dual of the cyclotomic $p$-units}

We shall actually be interested in the composition of $\Psi_K$ with
a map to a slightly different module arising from cyclotomic $S$-units.  We remark that 
$$
	\mf{X}_K \cong 
	\Hom(H^1(G_{K,S},\mu_{p^{\infty}}),\mu_{p^{\infty}}).
$$
Let $\mc{O}_{K,S}$ denote the ring of $S$-integers of $K$.
The direct limit of the sequences \eqref{kummer1} over $r$ and $E$
provides an exact sequence
$$
	0 \to \lim_{\rightarrow} \mc{O}_{K,S}^{\times}/\mc{O}_{K,S}^{\times p^r} \to
	H^1(G_{K,S},\mu_{p^{\infty}}) \to A_{K,S} \to 0
$$
for every $n \ge 1$,
where 
\label{AKS}
$A_{K,S}$ denotes the direct limit of the $p$-parts of the $S$-class groups
of number fields in $K$.  As before, let 
$\mc{E}_K$ denote the pro-$p$ completion of
$\mc{O}_{K,S}^{\times}$.  Note that
$$
	\Hom_{\cts}(\mc{E}_K,\zp) \cong 
	\Hom(\mc{O}_{K,S}^{\times},\zp) \cong \Hom(\mc{O}_{K,S}^{\times}
	\otimes_{\Z} \qp/\zp,\qp/\zp).
$$
Thus, we obtain an exact sequence
\begin{equation} \label{kummercons}
	0 \to \Hom(A_{K,S},\mu_{p^{\infty}}) \to \mf{X}_K \to
	\Hom_{\cts}(\mc{E}_K,\zp(1)) \to 0.
\end{equation}

Now consider the pro-$p$ completion
\label{CK}
$\mc{C}_K$ of the group 
of cyclotomic $S$-units in $K$, which is to say the pro-$p$ completion of the group generated by the $1-\xi$ with $\xi \in \mu_{Np^{\infty}}$, $\xi \neq 1$.  Let us set 
\label{mfCK}
$\mf{C}_K = \Hom_{\cts}(\mc{C}_K,\zp(1))$.  
By \eqref{kummercons}, we have a homomorphism
$$
	q \colon \mf{X}_K \to \mf{C}_K.
$$

Note that we may decompose a $\Lambda_N$-module $A$ into its $(\pm 1)$-eigenspaces
\label{Apm3} 
$A^{\pm}$ for the action of $-1 \in (\Z/Np\Z)^{\times}$.

\begin{remark}
	If $N = 1$ (or $2$), the map $q^- \colon \mf{X}_K^- \to \mf{C}_K^-$ 
	is an isomorphism if and only if Vandiver's conjecture holds.
\end{remark}

\begin{proposition} \label{CKprop}
	There exists a injection
	$$
\label{ThetaK}
		\Theta_K \colon \mf{C}_K^- \hookrightarrow (\Lambda_N^{\star})^-
	$$
	of $\Lambda_N$-modules, canonical up to the choice of $\iota$,
	such that
	$$
		\Theta_K(\phi) \otimes \zeta =
		\lim_{\substack{\leftarrow\\r}} \sum_{i=1}^{Np^r-1}
		\phi(1-\zeta_{Np^r}^i) \otimes [i]_r
	$$
	for all $\phi \in \mf{C}_K^-$.
\end{proposition}

\begin{proof}
	For $r \ge 1$, 
	let $C_r$ denote the $p$-completion of the cyclotomic $S$-units of
	$F_r = F(\mu_{p^r})$. 
	There are obvious surjections
	$$
		\psi_r \colon \zp[(\Z/Np^r\Z)^{\star}] \to C_r
	$$
	given by $\psi_r([i]_r) = 1-\zeta_{Np^r}^i$ for $i$ not divisible
	by $Np^r$.  
	These are compatible in the sense that
	$$
		\psi_r([i]_r) = \prod_{k=0}^{p^{s-r}-1}	\psi_s([i+kNp^r]_s)
	$$
	for any $s \ge r$.
	Note that we have
	$$
		\Hom_{\zp}(\zp[(\Z/Np^r\Z)^{\star}],\zp) \cong \zp[(\Z/Np^r\Z)^{\star}]
	$$
	via 
	$$
		\varphi \mapsto \sum_{i = 1}^{Np^r-1} \varphi([i]_r)[i]_r,
	$$ 
	and these are compatible in the
	sense that 
	$$
		\sum_{i=1}^{Np^r-1} \varphi_r([i]_r)[i]_r
		= \sum_{i=1}^{Np^{r+1}-1} \varphi_{r+1}([i]_{r+1})[i]_r
	$$
	if 
	$$
		\varphi_r \in \Hom_{\zp}(\zp[(\Z/Np^r\Z)^{\star}],\zp)
	$$ 
	for
	each $r \ge 1$ satisfy
	$$
		\varphi_r([i]_r) = \sum_{k=0}^{p^{s-r}-1}		
		\varphi_s([i+kNp^r]_s)
	$$
	for all $i$ not divisible by $Np^r$.  Hence, we have injections
	$$
		\psi_r^* \colon \Hom_{\zp}(C_r,\zp) 
		\hookrightarrow \zp[(\Z/Np^r\Z)^{\star}]
	$$
	dual to the $\psi_r$ that are compatible under restriction on the left
	and projection on the right.
	We have
	$$
		\mf{C}_K \cong \lim_{\substack{\leftarrow\\r}}\, \Hom_{\zp}(C_r,\zp(1)),
	$$
	where the inverse limit is taken with respect to restriction maps,
	so our injection $\Theta_K$ can be taken to be the restriction to
	$\mf{C}_K^-$ of the inverse limit of the $\psi_r^* \otimes 
	\zeta^{\otimes -1}$.
\end{proof}

\begin{remark} 
	In fact, $\mf{C}_K^-$ is isomorphic to $\Lambda_N^-$, but
	it is the injection of Proposition \ref{CKprop} that is most natural in our
	setting.  
	The map $\Theta_K$ is an isomorphism in its $\psi$-eigenspace 
	for any odd, primitive character $\psi$ of $(\Z/Np\Z)^{\times}$.
\end{remark}

Let 
\label{phi2}
$\phi_2$ be the composition of $q^- \colon \mf{X}_K^- \to \mf{C}_K^-$ with $\Theta_K$.
	
\subsection{The reciprocity map and the $L$-function} \label{recLfn}

Let $M$ be a positive divisor of $N$.
By Corollary \ref{LN-}, the image of $\mc{L}_{N,M}^{\star}$
in $\mc{Z}_N \cotimes \Lambda_N^{\star}$ is actually contained in
$\mc{Y}_N \cotimes \Lambda_N^{\star}$.
Let
$$
\label{overlineLNMstar}
	\overline{\mc{L}_{N,M}^{\star}} \in 
	\mc{Y}_N^-/\I \mc{Y}_N^- \otimes_{\zp} 
	(\Lambda_N^{\star})^-
$$  
denote the projection of $\mc{L}_{N,M}^{\star}$ to the latter Galois module.
(Recall that the Galois modules $\mc{H}_1^{\et}(N)$ and $\mc{H}^1_{\et}(N)(1)$ are canonically isomorphic, so $\overline{\mc{L}_{N,M}^{\star}}$ depends upon our fixed choice of $\iota$.)
We denote by $\Xi_N$ the homomorphism of
$\Lambda_N \cotimes \Lambda_N$-modules
$$
\label{XiN}
	\Xi_N = \phi_1 \otimes \phi_2
	\colon X_K^{\circ} \otimes_{\zp} 
	\mf{X}_K^- \to \mc{Y}_N^-/\I\mc{Y}_N^- \otimes_{\zp} 
	(\Lambda_N^{\star})^-
$$
resulting from Propositions \ref{1stisom} and \ref{CKprop}.

Recall that $H^2_S(K,\zp(1))^{\circ} \cong X_K^{\circ}$.
We will use 
\label{PsiKcirc}
$\Psi_K^{\circ}$ to denote the projection of $\Psi_K$ (see Section \ref{recmap}) to a map
$$
	\Psi_K^{\circ} \colon \mc{U}_K \to
	X_K^{\circ} \otimes_{\zp} \mf{X}_K^-.
$$
Again, let $1-\zeta^M \in \mc{U}_K$ denote the norm compatible sequence
$(1-\zeta_{Np^r}^M)_r$ of $S$-units.  We are now ready to state our main conjecture.

\begin{conjecture} \label{mainconj}
	We have 
	$$
		\Xi_N(\Psi_K^{\circ}(1-\zeta^M)) = \overline{\mc{L}_{N,M}^{\star}}. 
	$$
\end{conjecture}

In fact, Conjecture \ref{mainconj} is equivalent to our earlier conjectures
relating cup products and modular symbols.

\begin{proposition}
	Conjectures \ref{symbolconj} and \ref{mainconj} are equivalent.
\end{proposition}

\begin{proof}
	Let $Q \ge 2$ be a divisor of $Np^r$ for some $r \ge 0$, and let
	$i \in \Z$ with $Q \nmid i$.
	Define
	$$
		\pi_{i,Q} \colon \Lambda_N^{\star} \to \zp(1)
	$$
	by
	$$
		\pi_{i,Q}([j]) = \begin{cases} \zeta & j \equiv i \bmod Q \\
		1 & j \not\equiv i \bmod Q, \end{cases}
	$$
	which induces a map on $(\Lambda_N^{\star})^-$, viewing it
	as a submodule.
		
	We claim that
	$$
		\pi_{i,Q} \circ \phi_2 = \pi_{1-\zeta_Q^i}
	$$
	on $\mf{X}_K^-$, where $\pi_{1-\zeta_Q^i}$ is as in Section \ref{recmap}.
	Let $\sigma \in \mf{X}_K^-$.  
	Recall that $\phi_2$ is the composite of the map 
	$\mf{X}_K^- \to \mf{C}_K^-$ with $\Theta_K$.
	Then, by Proposition \ref{CKprop}, we have
	\begin{align*}
		\pi_{i,Q} \circ \phi_2(\sigma) 
		&= \pi_{i,Q}
		\Bigl(\lim_{\substack{\leftarrow\\s}} \sum_{j=1}^{Np^s-1}
		(\pi_{1-\zeta_{Np^s}^j}(\sigma) \otimes \zeta^{\otimes -1})[j]_s\Bigr)\\
		&= \lim_{\substack{\leftarrow\\s}}
		\prod_{\substack{j=1\\j \equiv i \bmod Q}}^{Np^s-1} 
		\pi_{1-\zeta_{Np^s}^j}(\sigma)\\
		&= \pi_{1-\zeta_Q^i}(\sigma)
	\end{align*}
	for $\sigma \in \mf{X}_K^-$, as desired.

	For any positive divisor $M$ of $N$, we have by definition that
	$$
		(1 \otimes \pi_{1-\zeta_Q^i})(\Psi_K^{\circ}(1-\zeta^M)) 
		= (1-\zeta_Q^i,1-\zeta^M)_{K,S}^{\circ}.
	$$
	It follows that
	\begin{equation} \label{Lident}
		(1 \otimes \pi_{i,Q})(\Xi_N(\Psi_K^{\circ}(1-\zeta^M)))
		= \Upsilon_K((1-\zeta_Q^i,1-\zeta^M)_{K,S}^{\circ}).
	\end{equation}
	
	Assume now that $(Np^r/Q \cdot i, M) = 1$.
	Since $U_p-1 \in \I$, we have 
	\begin{align*}
		(1 \otimes \pi_{i,Q})(\overline{\mc{L}_{N,M}^{\star}}) &= 
		\lim_{\substack{\leftarrow\\r}} 
		\sum_{\substack{j=1\\(j,M) = 1}}^{Np^r-1} 
		\bar{\xi}_r(j:M) \otimes \zeta^{\otimes -1} 
		\otimes \pi_{i,Q}([j]) \\
		&= \lim_{\substack{\leftarrow\\r}} 
		\sum_{k=0}^{(Np^r/Q)-1} \bar{\xi}_r(i+kQ:M).
	\end{align*}
	As in \eqref{maninredsum}, we have	
	$$
		\bar{\xi}_r((Np^r/Q)i:M)
		= \sum_{k=0}^{(Np^r/Q)-1} \bar{\xi}_r(i+kQ:M).
	$$
	Hence, we have that
	\begin{equation} \label{Rident}
		(1 \otimes \pi_{i,Q})(\overline{\mc{L}_{N,M}^{\star}}) = 
		\bar{\xi}((N/Q)i:M).
	\end{equation}
	
	Putting \eqref{Lident} and \eqref{Rident} together, 
	Conjecture \ref{mainconj} yields
	$$
		\Upsilon_K((1-\zeta_Q^i,1-\zeta^M)_{K,S}^{\circ})
		= \bar{\xi}((N/Q)i:M)
	$$
	for all $i$ and $Q$ as above.
	As any symbol $(1-\zeta_{Np^s}^u,1-\zeta^v)_{K,S}^{\circ}$ as in Conjecture
	\ref{symbolconj} is a Galois conjugate of one of the above form,
	and noting \eqref{maninrels4} and the fact that 
	$$
		\Upsilon_K \circ \sigma_j = \langle j \rangle^{-1} \circ \Upsilon_K,
	$$ 
	Conjecture \ref{mainconj} implies Conjecture \ref{symbolconj}.
	 
	Conversely, fix $M$ and take $i$ and $Q$ as before with $(Np^r/Q \cdot i,M)
	=1$.  Conjecture \ref{symbolconj} along with \eqref{Lident} 
	and \eqref{Rident} then imply that
	\begin{equation} \label{onedir}
	 	(1 \otimes \pi_{i,Q})(\Xi_N(\Psi_K^{\circ}(1-\zeta^M)))
		= (1 \otimes \pi_{i,Q})(\overline{\mc{L}_{N,M}^{\star}}).
	\end{equation}
	If instead $(i,M) \neq 1$, then $\pi_{i,Q}$
	is trivial on symbols of the form $[j]$ with $j \in \Z_{p,N}^{\star}$ 
	coprime to $M$, and 
	\eqref{onedir} still holds with both sides being trivial.  
	Since the $\pi_{i,Q}$ with $Q = Np^r$ for some $r \ge 1$ 
	and $i \in \Z$ not divisible by $Q$
	generate $\Hom_{\cts}(\Lambda_N^{\star},\zp(1))$ topologically,
	we have the reverse implication as well.
\end{proof}

\section{Comparison with $p$-adic $L$-values} 

\subsection{Characters and cyclotomic units} \label{cyclunits}

From now on, we work with multiple characters of $\Z_{p,N}^{\times}$ at
once, so it is easiest to extend scalars to the ring
\label{ON}
$\mc{O}_N = \zp[\mu_{\varphi(N)p^{\infty}}]$.
Similarly to the appendix to \cite{me-paireis} (but with a larger ring), for a $\zp[(\Z/Np\Z)^{\times}]$-module
$A$ and character $\chi$ of $(\Z/Np\Z)^{\times}$, we define
$$
\label{Achi}
	A^{\chi} \cong A^{(\chi)} \otimes_{R_{\chi}} \mc{O}_N
$$
to be the $\chi$-eigenspace of $A \otimes_{\zp} \mc{O}_N$ as an 
$\mc{O}_N[(\Z/Np\Z)^{\times}]$-module.
Furthermore, for any homomorphism $\alpha \colon A \to B$ of $\zp[(\Z/Np\Z)^{\times}]$-modules, we have an induced map
$$
	\alpha^{\chi} \colon A^{\chi} \to B^{\chi},
$$ 
which we may also view as a map from $A \otimes_{\zp} \mc{O}_N$ to $B^{\chi}$ factoring through
$A^{\chi}$.
Let $\epsilon_{\chi} \colon A \to A^{\chi}$ denote the idempotent
$$
\label{epsilonchi}
	\epsilon_{\chi} = \frac{1}{\varphi(Np)}\sum_{i \in (\Z/Np\Z)^{\times}}
	\chi(i)^{-1}[i]_1 \in \mc{O}_N[(\Z/Np\Z)^{\times}].
$$
For notational purposes, we extend these definitions to any function $\chi$ of $\Z_{p,N}^{\times}$ by using the restriction of $\chi$ to $(\Z/Np\Z)^{\times}$.

We extend 
\label{kappa2}
$\kappa$ multiplicatively to $\Z_{p,N}$ by setting $\kappa(p) = p$ and, if $l$ is a prime dividing $N$, taking $\kappa(l)$ to be the value one obtains by viewing $\kappa$ as a character on $\zp^{\times}$, which contains $l$.
Suppose now that $\chi \colon \Z_{p,N} \to \qpbar$ has the form $\chi = \psi\kappa^t$ for some $t \ge 0$, where $\psi$ arises as the continuous extension of a not necessarily primitive Dirichlet character on $\Z$ of period dividing $Np^r$ for some $r \ge 1$.
Following \cite{gs}, we refer to such a character as an {\em arithmetic character}.
	For such a $\chi$, we let \label{fchi} $f_{\chi}$ denote the prime-to-$p$ part of the 
	period of the restriction of $\psi$ to $\Z$.
	We consider 
	\label{omega2}
	$\omega$ as an arithmetic character by taking its unique extension with
	$f_{\omega} = 1$.

Let $\psi$ be a finite, even arithmetic character on $\Z_{p,N}$.
Fix $t \ge 1$ and consider positive integers $M$ dividing $Np$ and $Q$ dividing $N$.  Consider the products
\label{etaMrt}
\begin{align*}
	\eta_{M,r,t}^{\psi} &= 
	\prod_{\substack{i=1\\(i,M)=1}}^{Np^r-1}
	(1-\zeta_{Np^r}^i)^{\psi\kappa^{t-1}(i)} \\
	\alpha_{r,t}^{Q,\psi} &= \prod_{\substack{i=1\\(i,Np)=1}}^{Np^r-1}
	(1-\zeta_{Qp^r}^i)^{\psi\kappa^{t-1}(i)}
\end{align*}
for $r \ge 1$.  (Note the abuse of notation here: these elements lie
in $\mc{C}_K \otimes_{\zp} \mc{O}_N$, and so we allow exponents in $\mc{O}_N$.)
In fact, we may consider $\alpha_{r,t}^{Q,\psi}$ with
the same definition for any $t \in \zp$.
The $\eta_{M,r,t}^{\psi}$ satisfy
$$
	\eta_{M,r+1,t}^{\psi}(\eta_{M,r,t}^{\psi})^{-1} \in \mc{C}_K^{p^r}
$$
(for sufficiently large $r$) and similarly for the $\alpha_{r,t}^{Q,\psi}$.
Let us consider the limits
\begin{eqnarray} \label{unitdef}
	\eta_{M,t}^{\psi} = \lim_{r \to \infty} \eta_{M,r,t}^{\psi} & \text{and} &
	\alpha_t^{Q,\psi} = \lim_{r \to \infty} \alpha_{r,t}^{Q,\psi}.
\end{eqnarray}
\label{etaMt}
The Galois automorphism $\sigma_j$ corresponding to 
$j \in \Z_{p,N}^{\times}$ satisfies
$$
	\sigma_j(\eta_{M,t}^{\psi}) = 
	(\eta_{M,t}^{\psi})^{\psi^{-1}\kappa^{1-t}(j)},
$$
and likewise for the $\alpha_t^{Q,\psi}$, so these are elements of $\mc{C}_K^{\psi^{-1}}$.
Note that $\alpha_t^{Q,\psi} = 1$ if the prime-to-$p$ part of the conductor of $\psi$
does not divide $Q$.
Finally, we set $\alpha_t^{\psi} = \alpha_{t}^{N,\psi}$  
\label{alphat}
and use corresponding notation with ``$r$'' in the subscript for the $r$th terms which have these limits.
The limit elements compare as follows.
\begin{lemma} \label{compareunits}
	We have the following equalities:
	\begin{enumerate}
	\item[a.]
	$
		\alpha_t^{\psi} = (\eta_{M,t}^{\psi})^{\prod_{l \mid Np,\, l \nmid M} 
		(1- \psi\kappa^{t-1}(l))},
	$ 
	and
	\item[b.]
	$
		\alpha_t^{\psi} = 
		(\alpha_t^{Q,\psi})^{\frac{\varphi(Q)}{\varphi(N)}\prod_{l \mid N,\, l \nmid Q} 
		(1-\psi\kappa^{t-1}(l))}
	$
	if $f_{\psi} \mid Q$.
	\end{enumerate}
	Here, the products are taken over primes $l$.
\end{lemma}

\begin{proof}
	Set $\chi = \psi\kappa^{t-1}$.  Let us also consider
	$$
		\beta_{M,r,t}^{D,\psi} =
		\prod_{\substack{i=1\\(i,M)=1}}^{Dp^r-1}
		(1-\zeta_{Dp^r}^i)^{\chi(i)}
	$$
	for any $D$ dividing $N$ and $M$ dividing $Dp$, as well as the
	limits $\beta_{M,t}^{D,\psi}$ that exist when $f_{\psi} \mid D$.  	
	We claim that if $f_{\psi} \mid D$, then
	\begin{equation} \label{claimone}
		\beta_{M,t}^{D,\psi} = \eta_{M,t}^{\psi}.
	\end{equation}
	To see this, note that if the period of $\psi$ on $\Z$ divides $Dp^r$, then
	\begin{align*}
		\prod_{\substack{j=1\\(j,M)=1}}^{Dp^r-1}
		(1-\zeta_{Dp^r}^{j})^{\chi(j)}
		&= \prod_{\substack{j = 1\\(j,M) =1}}^{Dp^r-1}\prod_{k=0}^{N/D-1}
		(1-\zeta_{Np^r}^{j+kDp^r})^{\chi(j)}\\
		&\equiv \prod_{\substack{j = 1\\(j,M) =1}}^{Dp^r-1}\prod_{k=0}^{N/D-1}
		(1-\zeta_{Np^r}^{j+kDp^r})^{\chi(j+kDp^r)} \bmod \mc{C}_K^{p^r}.
	\end{align*}
	Since $j+kDp^r$ is prime to $M$ if $j$ is, the latter term is 
	$\eta_{M,r,t}^{\psi}$.  The claim then follows by taking limits.

	Consider the formal identity
	\begin{equation} \label{formal}
		\sum_{\substack{d \mid Np,\, (d,M) = 1\\d \ge 1}} \mu(d)
		\sum_{\substack{i=1\\(i,M) = 1}}^{Np^r/d-1} [di]_r = 
		\sum_{\substack{i=1\\(i,Np) = 1}}^{Np^r-1} [i]_r,
	\end{equation}
	where $\mu$ denotes the M\"obius function.	
	It follows from our definitions that
	\begin{equation} \label{identone}
		\alpha_{r,t}^{\psi} = 
		\prod_{\substack{d \mid Np,\, (d,M) = 1\\d \ge 1}}
		(\beta_{M,r-\delta_d,t}^{N/g_d,\psi})^{\mu(d)\chi(d)},
	\end{equation}
	where $g_d = (d,N)$ and $\delta_d = \log_p(d/g_d)$.
	Note that $\chi(d) \neq 0$ implies that $f_{\psi} \mid (N/g_d)$.  
	Taking limits and applying \eqref{claimone}, we get
	$$
		\alpha_t^{\psi} = 
		(\eta_{M,t}^{\psi})^{	
		\sum_{\substack{d \mid Np,\, (d,M) = 1}} \mu(d)\chi(d)}.
	$$
	Applying
	\begin{equation} \label{muident}
		\prod_{\substack{l \mid Np,\, l \nmid M\\ l \text{ prime}}} (1-\chi(l))
		= \sum_{\substack{d \mid Np,\, (d,M) = 1\\d \ge 1}} \mu(d)\chi(d),
	\end{equation}
	we have part a.

	For part b, we assume that $f_{\psi} \mid Q$. Then, we have
	\begin{equation} \label{identtwo}
		\alpha_t^{\psi} 
		= (\eta_{Qp,t}^{\psi})^{\prod_{l \mid N,\,l \nmid Q} (1-\chi(l))}
		= (\beta_{Qp,t}^{Q,\psi})^{\prod_{l \mid N,\,l \nmid Q} 
		(1-\chi(l))}.
	\end{equation}
	using part a in the first step and \eqref{claimone} in the second.
	On the other hand, we have
	$$
		\alpha_{r,t}^{Q,\psi} \equiv 
		(\beta_{Qp,r,t}^{Q,\psi})^{\frac{\varphi(N)}{\varphi(Q)}}
		\bmod{\mc{C}_K^{p^r}}.
	$$
	Taking limits and combining this with \eqref{identtwo}, we obtain
	part b.
\end{proof}

\subsection{Special values} \label{specval}

Let $A$ be a finitely generated $\Lambda_N$-module.  Then for any 
arithmetic character $\chi$ on $\Z_{p,N}$ we have specialization maps
\label{tildechi3}
\begin{eqnarray*}
	\tilde{\chi} \colon A \cotimes \Lambda_N^{\star} \to A \otimes_{\zp}
	\mc{O}_N, &&
	a \otimes [j] \mapsto a \otimes \chi(j).
\end{eqnarray*}
For later use, we remark that an element of $A \cotimes \Lambda_N^{\star}$ is uniquely determined by its specializations. 

\begin{lemma} \label{specialize}
	Suppose that $A$ is $\zp$-torsion free.
	An element $a \in A \cotimes \Lambda_N^{\star}$ satisfies 
	$\tilde{\chi}(a) = 0$ for all (finite) arithmetic characters $\chi$ on $\Z_{p,N}$ if and only if $a = 0$.
\end{lemma}

\begin{proof}
	An element of $A \cotimes \Lambda_N^{\star}$ is nonzero if and only
	if, for every $\zp$-quotient $B$ of $A$ and $r \ge 1$,
	the image of $a$ in $B \otimes_{\zp} \zp[(\Z/Np^r\Z)^{\star}]$ is trivial.
	The problem then reduces to the case that $A = \zp$.  
	Now, choose $r \ge 1$ such that the image of $a$ in 
	$\zp[(\Z/Np^r\Z)^{\star}]$ is nonzero.  Since the primitive	
	Dirichlet characters
	of conductor dividing $Np^r$ form a basis of the space of $\qpbar$-valued
	congruence functions of period dividing $Np^r$, 
	there exists a primitive	
	Dirichlet character $\psi$ of conductor
	dividing $Np^r$ such that $\tilde{\psi}(a) \neq 0$.
\end{proof}

Let us compare $\mc{L}_{N,M}$ and $\mc{L}_{N,M}^{\star}$ for any positive
divisor $M$ of $N$.  

\begin{lemma} \label{compareL}
	For any arithmetic character $\chi$ on $\Z_{p,N}$ we have
	$$
		U_D\tilde{\chi}(\mc{L}_{N,M}) =
		\biggl(\prod_{\substack{l \mid Np,\, l \nmid M\\l \mr{\ prime}}} 
		(U_l-\chi(l))\biggr)
		\tilde{\chi}(\mc{L}_{N,M}^{\star}),
	$$
	where $D = \prod_{l \mid Np,\, l \nmid M} l$.
\end{lemma}

\begin{proof}
	Using \eqref{formal}, we have
	$$
		\sum_{\substack{j=1\\(j,Np)=1}}^{Np^r-1} \chi(j)U_p^{-r}\xi_r(j:M)
		= \sum_{\substack{d \mid Np\\(d,M)=1\\d \ge 1}}
		\mu(d)\sum_{\substack{j=1\\(j,M)=1}}^{Np^r/d-1} 
		\chi(dj)U_p^{-r}\xi_r(dj:M).
	$$
	If $\mu(d) \neq 0$ in the latter equation, then $d$ divides $D$, and we
	have
	$$
		U_D\xi_r(dj:M) = U_{D/d}\sum_{k=0}^{d-1}\xi_r(j+kNp^r/d:M)
	$$
	by \eqref{maninsum}.  Let $r$ be large enough that
	$$
		 \chi(dj) \equiv \chi(d)\chi(j + kNp^r/d) \bmod p^r\mc{O}_N.
	$$	
	Since $M \mid (Np^r/d)$, we then have
	$$
		U_D\sum_{\substack{j=1\\(j,Np)=1}}^{Np^r-1} \chi(j)U_p^{-r}\xi_r(j:M) 
		\equiv \sum_{\substack{d \mid Np\\(d,M)=1\\d \ge 1}} \mu(d)\chi(d)U_{D/d} 
		\sum_{\substack{j=1\\(j,M)=1}}^{Np^r-1} \chi(j)U_p^{-r}\xi_r(j:M),
	$$
	where the congruence is taken modulo $p^rs_r(H_1(X_1^r(N),C_1^r(N);\zp))^{\ord}$.
	We obtain
	$$
		U_D\tilde{\chi}(\mc{L}_{N,M}) =
		\biggl(\sum_{\substack{d \mid Np,\, (d,M)=1\\d \ge 1}} 
		\mu(d)\chi(d)U_{D/d}\biggr)
		\tilde{\chi}(\mc{L}_{N,M}^{\star})
	$$
	in the limit.
	By an analogous equation to \eqref{muident}, the result follows.
\end{proof}

For any $\chi$ on $\Z_{p,N}^{\times}$, let 
$$
\label{hchi}
	\mf{h}_{\chi} = (\mf{h} \otimes_{\zp} \mc{O}_N)/(\langle a \rangle -
	\chi\kappa^{-2}\omega^{-2}(a) \mid a \in \Z_{p,N}^{\times}).
$$
Then, for any $\mf{h}$-module $Z$, set 
\label{Zchi}
$Z_{\chi} = Z \otimes_{\mf{h}} 
\mf{h}_{\chi}$,
and let $P_{\chi} \colon Z \to Z_{\chi}$
\label{Pchi}
be the natural map.
Now, if $\alpha$ and $\chi$ are finite 
	arithmetic
characters on $\Z_{p,N}$ 
and $k, s \in \zp$, then we define
$$
\label{LpM}
	L_{p,M}(\xi,\alpha,k,\chi,s) = 
	P_{\alpha\kappa^k}(\widetilde{\chi\kappa^{s-1}}
	\left(\mc{L}_{N,M})\right) \in H_1(N)_{\alpha\kappa^k}.
$$
If $s$ is a positive integer, then we set
$$
	L_{p,M}^{\star}(\xi,\alpha,k,\chi,s) = 
	P_{\alpha\kappa^k}(\widetilde{\chi\kappa^{s-1}}\left(\mc{L}_{N,M}^{\star})
	\right) \in \mc{H}_1(N)_{\alpha\kappa^k}.
$$
It follows from \eqref{maninrels3} and \eqref{maninrels4} that $L_{p,M}^{\star}(\xi,\alpha,k,\chi,s)$		
is zero if $\alpha$ is odd.  Lemma \ref{compareL} has the following immediate corollary.

\begin{corollary} \label{compareLp}
	Let $\alpha$ and $\chi$ be finite 
		arithmetic	
	characters on
	$\Z_{p,N}$, let $k \in \zp$, and let $s$ be a positive integer.
	Then we have
	$$
		U_D L_{p,M}(\xi,\alpha,k,\chi,s) 
		= \biggl(\prod_{\substack{l \mid Np,\, l \nmid M\\l\,\mr{prime}}} 
		(U_l-\chi\kappa^{s-1}(l))\biggr)
		L_{p,M}^{\star}(\xi,\alpha,k,\chi,s),
	$$
	with $D$ as in Lemma \ref{compareL}.
\end{corollary}

Let us abbreviate the standard $p$-adic $L$-function $L_{p,1}$ by $L_p$.  
\label{Lp}
This has the following functional equation:

\begin{lemma} \label{fnlequation}
	We have
	$$
		L_p(\xi,\alpha,k,\chi,s) = -\chi(-1)L_p(\xi,\alpha,k,\alpha\chi^{-1}\omega^{-2},k-s).
	$$
\end{lemma}

\begin{proof}
	We may assume that $\alpha$ is even.
	The result then follows directly from \eqref{maninrels1} and \eqref{maninrels4}, which
	yield the identity
	\begin{equation} \label{diamident}
		P_{\alpha\kappa^k}(\xi_r(j:1)) 
		= -\alpha^{-1}\omega^2\kappa^{-k+2}(j)P_{\alpha\kappa^k}(\xi_r(-j^{-1}:1))
	\end{equation}
	for $j$ prime to $Np$.
	This in turn implies
	$$
		P_{\alpha\kappa^k}\biggl(\sum_{\substack{j=1\\(j,Np)=1}}^{Np^r-1}
		\chi\kappa^{s-1}(j)\xi_r(j:1)\biggr)
		= -\chi(-1)P_{\alpha\kappa^k}\biggl(\sum_{\substack{j=1\\(j,Np)=1}}^{Np^r-1}
		\chi\alpha^{-1}\omega^2\kappa^{s-k+1}(j)\xi_r(j^{-1}:1)\biggr),
	$$
	yielding the desired result in the inverse limit.
\end{proof}

We also have the following.

\begin{lemma} \label{compareLM}
	Let $Q = N/M$, and suppose that $f_{\alpha\chi^{-1}} \mid Q$.  Then
	$$
		U_D L_p(\xi,\alpha,k,\chi,s)
		= \frac{\varphi(Q)}{\varphi(N)}U_M
		\biggl(\prod_{\substack{l \mid N\\l \nmid Q}}
		(U_l-\alpha\chi^{-1}\omega^{-2}\kappa^{k-s-1}(l))\biggr) 
		L_{p,M}(\xi,\alpha,k,\chi,s)
	$$ 
	for $D = \prod_{l \mid N,\, l \nmid Q} l$.
\end{lemma}

\begin{proof}
		First, let us remark that, while $\alpha\chi^{-1}$ is not a priori well-defined,
		we make it so by considering it as a finite arithmetic character by taking the extension
		of its restriction to $\Z_{p,N}^{\times}$ of minimal period.
	We compute the latter $L$-value.
	Let us define
	$$
		\mc{L}_N^{Qp} = \lim_{\leftarrow} 
		\sum_{\substack{j=1\\(j,Qp)=1}}^{Np^r-1} U_p^{-r}\xi_r(j:1) \otimes
		[j]_r.
	$$
	Just as in the proof of Lemma \ref{compareL}, we have
	$$	
		U_D\tilde{\beta}(\mc{L}_N)
		= \biggl(\prod_{\substack{l \mid N\\l \nmid Q}}
		(U_l-\beta(l))\biggr)\tilde{\beta}(\mc{L}_N^{Qp})
	$$
	for any arithmetic character $\beta$ on $\Z_{p,N}$. 
	Hence, setting 
	$$
		L_p^{Qp}(\xi,\alpha,k,\chi,s) = P_{\alpha\kappa^k}\left(
		\widetilde{\chi\kappa^{s-1}}(\mc{L}_N^{Qp})\right)
	$$ 
	in general, we obtain
	\begin{multline} \label{LpQpcomp}
		U_DL_p(\xi,\alpha,k,\alpha\chi^{-1}\omega^{-2},k-s) \\
		= \biggl(\prod_{\substack{l \mid N\\l \nmid Q}}
		(U_l-\alpha\chi^{-1}\omega^{-2}\kappa^{k-s-1}(l))\biggr)
		L_p^{Qp}(\xi,\alpha,k,\alpha\chi^{-1}\omega^{-2},k-s).
	\end{multline}
	Then, for $r$ sufficiently large,	
	\begin{align*}
		\sum_{\substack{j=1\\(j,Qp)=1}}^{Np^r-1}
		\alpha\chi^{-1}\omega^{-2}\kappa^{k-s-1}(j)\xi_r(j:1)
		&\equiv \sum_{\substack{j=1\\(j,Qp)=1}}^{Qp^r-1} 
		\alpha\chi^{-1}\omega^{-2}\kappa^{k-s-1}(j)
		\sum_{k=0}^{M-1} \xi_r(j+kQp^r:1)\\
		&= U_M\sum_{\substack{j=1\\(j,Qp)=1}}^{Qp^r-1} 
		\alpha\chi^{-1}\omega^{-2}\kappa^{k-s-1}(j) \xi_r(Mj:1)\\
		&\equiv \frac{\varphi(Q)}{\varphi(N)}U_M
		\sum_{\substack{j=1\\(j,Np)=1}}^{Np^r-1} 
		\alpha\chi^{-1}\omega^{-2}\kappa^{k-s-1}(j) \xi_r(Mj:1),
	\end{align*}
	where the congruences are taken modulo 
	$p^rs_r(H_1(X_1^r(N),C_1^r(N);\zp))^{\ord}$
	and we have used \eqref{maninsum} in the second step.
	Since $\xi_r(Mj:1) = \langle j \rangle^{-1}\xi_r(M:j^{-1})$, we 
	therefore have
	\begin{equation} \label{usefnl}
		L_p^{Qp}(\xi,\alpha,k,\alpha\chi^{-1}\omega^{-2},k-s)
		= -\chi(-1)	\frac{\varphi(Q)}{\varphi(N)}U_M
		L_{p,M}(\xi,\alpha,k,\chi,s),
	\end{equation}
	as desired.  The result now arises from \eqref{LpQpcomp} 
	by applying Lemma \ref{fnlequation}
	to its left-hand side and plugging in the
	result of \eqref{usefnl} on its right hand-side.
\end{proof}

\subsection{Cup products and special values}

Let $\psi$ and $\gamma$ be finite even
arithmetic characters on $\Z_{p,N}$, set $\theta = \omega\psi\gamma$,
and assume that $p \mid B_{1,\theta}$ and $\theta$ is primitive when restricted to $(\Z/Np\Z)^{\times}$.  
The pairing $(\,\cdot\,,\,\cdot\,)_{K,S}$ of Section \ref{Galois} induces an $\mc{O}_N$-bilinear map
$$
\label{cupKStp}
	(\,\cdot\,,\,\cdot\,)_{K,S}^{(\psi,\gamma)} \colon
	\mc{C}_K^{\psi^{-1}} \times \mc{U}_K^{\gamma^{-1}} \to 
	X_K^{\theta^{-1}}(1),
$$
where one should note that the inverses of the characters in question are well-defined on 
$\Z_{p,N}^{\times}$. Any element $b \in \mc{C}_K^{\psi^{-1}}$ induces a homomorphism 
$$
	\pi_b^{\psi} \colon \mf{X}_K \otimes_{\zp}	
	\mc{O}_N \to \mc{O}_N(1)
$$ 
factoring through $\mf{X}_K^{\omega\psi}$.  
By \cite[Lemma A.1]{me-paireis}, the
map	$(\Psi_K^{\circ})^{\gamma^{-1}}$ takes its image in 
$$
	(X_K^{\circ} \otimes_{\zp} \mf{X}_K^-)^{\gamma^{-1}} \cong 
	\bigoplus_{\substack{\chi \text{ even}\\ (\chi\theta^{-1}) \in \Sigma_{Np}}} 
	(X_K^{\chi\theta^{-1}} \otimes_{\mc{O}_N} \mf{X}_K^{\omega\psi\chi^{-1}}).
$$
It follows from \eqref{cuppsi} that
$$
	(b,u)_{K,S}^{(\psi,\gamma)} = (1 \otimes \pi_b^{\psi})((\Psi_K^{\circ})^{\gamma^{-1}}(u))
$$
for any $u \in \mc{U}_K^{\gamma^{-1}}$.

Finally, for any positive divisor $M$ of $N$, let 
$$
	(\overline{\mc{L}_{N,M}^{\star}})^{\ltheta} \in 
	(\mc{Y}_N^-/\I\mc{Y}_N^-)^{\ltheta}(1) \otimes_{\zp} 
	(\Lambda_N^{\star})^-
$$ 
denote the image of $\mc{L}_{N,M}^{\star}$ in this module.
(In what follows, we will view $\widetilde{\chi}((\overline{\mc{L}_{N,M}^{\star}})^{\ltheta})$
for an arithmetic character $\chi$ on $\Z_{p,N}$ as its image in 
$(\mc{Y}_N^-/\I\mc{Y}_N^-)^{\ltheta} \otimes_{R_{\theta}} \mc{O}_N$.)

\begin{proposition} \label{intermedconj}
	Conjecture \ref{symbolconj} is equivalent to the statement that	
	$$
	\Upsilon_K^{\omega\theta^{-1}}
	((\eta_{M,t}^{\psi},\epsilon_{\omega\psi\theta^{-1}}(1-\zeta^M))_{K,S}^{(\psi,\theta\omega^{-1}\psi^{-1})})
	= \widetilde{\psi\kappa^{t-1}}((\overline{\mc{L}_{N,M}^{\star}})^{\ltheta})
	$$
	for any positive integer $M$ dividing $N$, finite even arithmetic character $\psi$ on $\Z_{p,N}$,
	primitive odd character $\theta$ on $(\Z/Np\Z)^{\times}$ with $p \mid B_{1,\theta}$, and $t \ge 1$.
\end{proposition}

\begin{proof}
	For $u \in \Z[\frac{1}{p}]$ and $v \in \Z$ prime to $p$ with
	$(u,v,N)\Z[\frac{1}{p}] = \Z[\frac{1}{p}]$,
	let $\bar{\xi}^{\ltheta}(u:v)$ denote the projection of
	$\bar{\xi}(u:v)$ to $(\mc{Y}_N^-/\I\mc{Y}_N^-)^{\ltheta}$.
	We have
	\begin{equation} \label{writespec}
		\Upsilon_K^{\omega\theta^{-1}}((\eta_{M,r,t}^{\psi},
		\epsilon_{\omega\psi\theta^{-1}}(1-\zeta^M))_{K,S}^{(\psi,\theta\omega^{-1}\psi^{-1})})
		= \sum_{\substack{i=1\\(i,M) = 1}}^{Np^r-1} \psi\kappa^{t-1}(i)\Upsilon_K^{\omega\theta^{-1}}
		((1-\zeta_{Np^r}^i,1-\zeta^M)_{K,S}).
	\end{equation}
	Using Conjecture \ref{symbolconj}, this becomes
	$$
		\sum_{\substack{i=1\\(i,M) = 1}}^{Np^r-1} 
		\psi\kappa^{t-1}(i)\bar{\xi}^{\ltheta}(p^{-r}i:M)
		= \widetilde{\psi\kappa^{t-1}}\Biggl(\sum_{\substack{i=1\\(i,M)=1}}^{Np^r-1}
		\bar{\xi}^{\ltheta}(p^{-r}i:M) \otimes [i]_r\Biggr),
	$$
	and the limit over $r$ of the latter sum is 
	$(\overline{\mc{L}_{N,M}^{\star}})^{\ltheta}$ by its definition, noting
	\eqref{LNMstar}.  
	
	As for the reverse implication, noting \eqref{writespec}, 
	we may apply Lemma \ref{specialize} to obtain
	$$
		\Upsilon_K^{\omega\theta^{-1}}(
		(1-\zeta_{Np^r}^i,1-\zeta^M)_{K,S})
		= \bar{\xi}^{\ltheta}(p^{-r}i:M)
	$$
	for all $i$ not divisible by $Np^r$ and	
	primitive odd characters $\theta$ on $(\Z/Np\Z)^{\times}$.
	Conjecture \ref{symbolconj} follows immediately from this.
\end{proof}

For any character $\chi$ on $\Z_{p,N}^{\times}$
and $\zp$-algebra $\mc{O}$ containing the values of $\chi$, 
\label{ONchi}
let $\mc{O}(\chi)$ be $\mc{O}$ as an $\mc{O}$-module,
endowed with a $\chi$-action of the Galois group $\Z_{p,N}^{\times}$.
For any $\Lambda_N \otimes_{\zp} \mc{O}$-module $A$, let 
$A(\chi) = A \otimes_{\mc{O}} \mc{O}(\chi)$.
\label{A'chi}
Alternatively, if $A$ is a $\Lambda_N$-module, we set $A'(\chi) = A \otimes_{\zp} \mc{O}_N(\chi)$.
For notational convenience, we set
$$
	H^i_{\cts}(G_{\Q,S},\mc{O}_N(\chi)) = H^i_{\cts}(G_{\Q,S},R(\chi)) \otimes_R \mc{O}_N,
$$ 
for $i \ge 0$, where $R$ is the $\zp$-algebra generated by the values of $\chi$.
Let $k$ and $t$ be elements of $\zp$.
As in Lemma \ref{weight2map}, 
corestriction provides an isomorphism
$$
	(X_K'
	(\kappa^{k-1}\theta))_{\Gal(K/\Q)} \xrightarrow{\sim}
	H^2_{\cts}(G_{\Q,S},\mc{O}_N(\kappa^k\omega\theta)).
$$
It also provides a homomorphism
$$
	(\mc{U}_K'
	(\kappa^{k-t-1}\gamma))_{\Gal(K/\Q)} \to
	H^1_{\cts}(G_{\Q,S},\mc{O}_N(\kappa^{k-t}\omega\gamma))
$$
under which $1-\zeta$ is mapped to $\alpha_{k-t}^{\gamma}$.
We have a cup product map
$$
	H^1_{\cts}(G_{\Q,S},\mc{O}_N(\kappa^t\omega	\psi)) 
	\otimes_{\mc{O}_N} H^1_{\cts}(G_{\Q,S},\mc{O}_N(\kappa^{k-t}\omega\gamma))
	\xrightarrow{\cup} H^2_{\cts}(G_{\Q,S},\mc{O}_N(\kappa^k\omega\theta)).
$$

Define $\mf{h}_{\chi}^*$ and $P_{\chi}^*$ analogously to $\mf{h}_{\chi}$ and $P_{\chi}$ using the adjoint diamond operators.
We then let 
\label{Ichi}
$I_{\chi}$ denote the image of $\I$ in $(\mf{h}^*_{\chi})_{\mf{m}}$ and set 
\label{Ychi}
$Y_{\chi} = P_{\chi}^*(\mc{Y}_N)$.  We let 
\label{overlineLpMstar}
$\overline{L_{p,M}^{\star}(\xi,\alpha,k,\chi,s)}$ denote the image of $L_{p,M}^{\star}(\xi,\alpha,k,\chi,s)$ in
$Y_{\alpha\kappa^k}/I_{\alpha\kappa^k}Y_{\alpha\kappa^k}$ for any allowable $\alpha$, $\chi$, $k$, and $s$,
and similarly for $L_{p,M}$ and $L_p$.  
\label{overlineLp}

Note that $\Upsilon_K$ induces a homomorphism
$$
\label{nuchi}
	\nu_{\kappa^k\omega\theta} 
	\colon H^2_{\cts}(G_{\Q,S},\mc{O}_N(\kappa^k\omega\theta))
	\to (Y_{\kappa^k\omega\theta}^-/
	I_{\kappa^k\omega\theta}Y_{\kappa^k\omega\theta}^-)(\kappa^{k-1}\theta),
$$
and recall the definitions of the limits of $S$-units $\eta_{M,t}^{\psi}$ and $\alpha_t^{Q,\psi}$ from \eqref{unitdef}.  We make the following conjecture.

\begin{conjecture} \label{cupprodconj}
	Let $M$ be a positive integer dividing $N$, and let $\psi$ and $\gamma$
	be finite even arithmetic characters on $\Z_{p,N}$.  Set $\theta = \omega\psi\gamma$,
	suppose that $p \mid B_{1,\theta}$ and $\theta|_{(\Z/Np\Z)^{\times}}$ is	
	primitive, and let $k \in \zp$ and $t \ge 1$.  Then, we have
	$$
		\nu_{\kappa^k\omega\theta}(\eta_{M,t}^{\psi} \cup 
		\alpha_{k-t}^{N/M,\gamma})
		= \overline{L_{p,M}^{\star}(\xi,\omega\theta,k,\psi,t)}.
	$$
\end{conjecture}

One sees easily that $\eta_{M,t}^{\psi} \cup \alpha_{k-t}^{N/M,\gamma}$ is the image of $(\eta_{M,t}^{\psi},\epsilon_{\gamma^{-1}}(1-\zeta^M))^{(\psi,\gamma)}$ under corestriction.
Therefore, Proposition \ref{intermedconj} implies the following:

\begin{proposition} \label{cupconjequiv}
	Conjecture \ref{symbolconj} is equivalent to
	Conjecture \ref{cupprodconj}.
\end{proposition}

In terms of standard $p$-adic $L$-values, we have the following 
possibly weaker conjecture.

\begin{conjecture}  \label{weakcupprodconj}
	Let $\theta$ and $\psi$ be finite arithmetic characters on $\Z_{p,N}$ with 
	$\theta$ odd and $\psi$ even, suppose that $p \mid B_{1,\theta}$ and
	$\theta|_{(\Z/Np\Z)^{\times}}$ is	 primitive, and let $k, t \in \zp$.  Then we have
	$$
		\nu_{\kappa^k\omega\theta}(\alpha_{t}^{\psi} \cup \alpha_{k-t}^{\theta\psi^{-1}\omega^{-1}})
		= \overline{L_{p}(\xi,\omega\theta,k,\psi,t)}.
	$$
\end{conjecture}

\begin{proposition}
	Conjecture \ref{cupprodconj} for $M = 1$ implies Conjecture \ref{weakcupprodconj}, 
	and the converse implication holds if $\mc{Y}_N^-/\I \mc{Y}_N^-$ is $p$-torsion free.
\end{proposition}

\begin{proof}
	We will prove slightly more than what is claimed.
	Suppose first that $t \ge 1$.  Pick $M$ such that $f_{\gamma} \mid (N/M)$, e.g., $M = 1$.
	By Lemma \ref{compareunits}, we have that 
	$\alpha_{t}^{\psi} \cup \alpha_{k-t}^{\gamma}$ equals
	$$
		\frac{\varphi(N/M)}{\varphi(N)}
		\biggl(\prod_{\substack{l \mid Np\\l \nmid M}} 
		(1-\psi\kappa^{t-1}(l)) \biggr)
		\biggl(\prod_{\substack{l \mid N\\
		l \nmid N/M}} (1-\gamma\kappa^{k-t-1}(l))\biggr)
		(\eta_{M,t}^{\psi} \cup \alpha_{k-t}^{N/M,\gamma}).
	$$
	On the other hand, by Corollary \ref{compareLp} and Lemma \ref{compareLM},
	we have that $\overline{L_{p}(\xi,\omega\theta,k,\psi,t)}$
	equals
	$$
		\frac{\varphi(N/M)}{\varphi(N)}
		\biggl(\prod_{\substack{l \mid Np\\l \nmid M}} 
		(1- \psi\kappa^{t-1}(l)) \biggr)
		\biggl(\prod_{\substack{l \mid N\\
		l \nmid N/M}} (1- \gamma\kappa^{k-t-1}(l))\biggr)
		\overline{L_{p,M}^{\star}(\xi,\omega\theta,k,\psi,t)}.
	$$
	Conjecture \ref{cupprodconj} then immediately implies Conjecture \ref{weakcupprodconj}
	for $t \ge 1$.  The general case then follows by taking limits using a sequence of positive integers 
	converging to $t$, since $\alpha_t^{\psi}$, $\alpha_{k-t}^{\gamma}$, and 
	$L_p(\xi,\omega\theta,k,\psi,t)$ vary continuously with $t$.
	
	Conversely, suppose that
	$1 - \gamma(l)$ lies in $\mc{O}_N^{\times}$ for all $l \mid N$ with 
	$l \nmid (N/M)$.  This occurs, of course, whenever $M = 1$.
	Conjecture \ref{weakcupprodconj} then implies that
	$$
		\nu_{\kappa^k\omega\theta}(\alpha_t^{\psi} \cup
		\alpha_{k-t}^{N/M,\gamma}) = \overline{L_{p,M}(\xi,\omega\theta,k,\psi,t)}
	$$
	for all $k \in \zp$ and $t \ge 1$ 
	(again by Lemmas \ref{compareunits} and \ref{compareLM}), and hence that
	$$
		\Upsilon_K^{\omega\theta^{-1}}((\alpha_t^{\psi},\epsilon_{\gamma^{-1}}
		(1-\zeta^M))_{K,S}^{(\psi,\gamma)})
		= \widetilde{\gamma\kappa^{t-1}}((\overline{\mc{L}_{N,M}})^{\ltheta}).
	$$
	Next, note that 
	$$
		(\alpha_t^{\psi},\epsilon_{\gamma^{-1}}
		(1-\zeta^M))_{K,S}^{(\psi,\gamma)}
		= \biggl(\prod_{\substack{l \mid Np\\l \nmid M}} 
		(1- \psi\kappa^{t-1}(l)) \biggr)(\eta_{M,t}^{\psi},\epsilon_{\gamma^{-1}}
		(1-\zeta^M))_{K,S}^{(\psi,\gamma)},
	$$
	and
	$$
		\widetilde{\psi\kappa^{t-1}}(\overline{\mc{L}_{N,M}})
		= \biggl(\prod_{\substack{l \mid Np\\l \nmid M}} 
		(1-\psi\kappa^{t-1}(l)) \biggr)\widetilde{\psi\kappa^{t-1}}(\overline{\mc{L}_{N,M}^{\star}}).
	$$
	We have that $\psi\kappa^{t-1}(l) \neq 1$ for $t \ge 2$
	and $l \mid Np$, since $\kappa^{t-1}(l)$ is in this case 
	a nontrivial element of $1+p\zp$ and $\psi(l)$ is either $0$ or 
	a root of unity in $\mc{O}_N$.  Since $X_K^{\theta^{-1}}$ and, by assumption, 
	$
	(\mc{Y}_N^-/\I\mc{Y}_N^-)^{\ltheta}
	$ 
	contain no $p$-torsion, it follows that
	$$
		\Upsilon_K^{\theta^{-1}}((\eta_{M,t}^{\psi},\epsilon_{\gamma^{-1}}(1-\zeta^M))_{K,S}^{(\psi,\gamma)})
		= \widetilde{\psi\kappa^{t-1}}((\overline{\mc{L}_{N,M}^{\star}})^{\ltheta})
	$$
	for all $t \ge 2$.  Conjecture \ref{cupprodconj} for our
	our particular $M$, $\psi$ and $\gamma$
	then follows for $t = 1$ as well by Lemma \ref{specialize}.
\end{proof}

\begin{remark}
	If $p \nmid B_{1,\theta^{-1}}$ for a given $\theta$ with $p \mid B_{1,\theta}$, then 
	Proposition \ref{1stisom} implies that $(\mc{Y}_N^-/\I\mc{Y}_N^-)^{\ltheta}$ is
	$p$-torsion free.  As in the remark at the end of Section \ref{mapback},
	our conjectures imply that $\mc{Y}_N^-/\I\mc{Y}_N^-$ is isomorphic to $X_K^{\circ}$,
	so we expect it to be $p$-torsion free in general.
\end{remark}

\twocolumn[\large {\bf Index of Notation}\vspace{1ex}]
\begin{ilist}
	\item Bernoulli numbers
	\begin{ilist2}
		\item $B_{1,\theta}$, \pr{B1theta}
	\end{ilist2}
	\item characters
	\begin{ilist2}
		\item $\omega$, $\kappa$, \sr{omega}{kappa}, \sr{kappa2}{omega2}
	\end{ilist2}
	\item comparison maps
	\begin{ilist2}
	    	\item $\phi_1'$, $\phi_1$, \sr{phi1'}{phi1}
   	    	\item $\Upsilon_K$, \pr{UpsilonK}
    		\item $\nu_r$, \pr{nur}
	    	\item $\varpi_r$, \pr{varpir}
	    	\item $\Xi_N$, \pr{XiN}
	    	\item $\nu_{\chi}$, \pr{nuchi}
	\end{ilist2}
	\item complex embedding
	\begin{ilist2}
	    	\item $\iota$, \pr{iota}
	\end{ilist2}
	\item conductors
	\begin{ilist2}
	    	\item $f_{\chi}$, \pr{fchi}
	\end{ilist2}
	\item conjugacy classes
	\begin{ilist2}   	 
		\item $(\chi)$, $\Sigma$, $\Sigma_{Np}$, \sr{Sigma}{SigmaNp}
	\end{ilist2}
	\item cup product pairings
	\begin{ilist2}   	 
	    	\item $(\,\cdot\,,\,\cdot\,)_{E,S}$, $(\,\cdot\,,\,\cdot\,)_{K,S}$, 
			\sr{cupES}{cupKS}
	    	\item $(\,\cdot\,,\,\cdot\,)_{K,S}^{\circ}$, \pr{cupKScirc}
	    	\item $(\,\cdot\,,\,\cdot\,)_{F_r,S}^{\circ}$, \pr{cupFr}
	    	\item $(\,\cdot\,,\,\cdot\,)_{K,S}^{(\psi,\gamma)}$,\pr{cupKStp}
	\end{ilist2}
	\item cyclotomic fields
	\begin{ilist2} 
	    	\item $F$, $K$, \sr{F}{K}
	   	\item $F_r$, \pr{Fr}
	\end{ilist2}
	\item eigenspaces and maps
	\begin{ilist2} 
	    	\item $A^{\pm}$, \pr{Apm}, \pr{Apm2}, \pr{Apm3}
	    	\item $A^{\circ}$, $A^{(\chi)}$, \sr{Acirc}{A(chi)}
	    	\item $\varepsilon$, \pr{varepsilon}
	    	\item $Z^{\ltheta}$, \pr{Ztheta}
	    	\item $A^{\chi}$, $\alpha^{\chi}$, $\epsilon_{\chi}$, \sr{Achi}{epsilonchi}
	    	\item $Z_{\chi}$, $P_{\chi}$, \sr{Zchi}{Pchi}
	\end{ilist2}
	\item Eisenstein ideals
	\begin{ilist2} 
	    	\item $\I$, $\mf{m}$, $\mf{I}$, $\mf{M}$, \sr{I}{M}
	   	 \item $I_r$, \pr{Ir}
	    	\item $I_{\chi}$, \pr{Ichi}
	\end{ilist2}
	\item Galois and related groups
	\begin{ilist2} 
	    \item $G_{E,S}$, $\mf{X}_E$, $X_{K,S}$, \sr{GES}{XKS}
	    \item $X_K$, \pr{XK}
	    \item $A_{K,S}$, $\mf{C}_K$, \sr{AKS}{mfCK}
	\end{ilist2}
	 \item Galois cohomology
	\begin{ilist2} 
	    	\item $H^i_S(K,\mc{T})$, \pr{HiS}
	    	\item $H^2_{\cts}(G_{F_r,S},\zp(2))^{\circ}$, \pr{H2circ}
	\end{ilist2}
	\item Hecke algebras
	\begin{ilist2} 
	    \item $\mf{h}_r$, $\mf{H}_r$, \pr{Hr}
	    \item $\mf{h}_r^{\ord}$, $\mf{H}_r^{\ord}$, $\mf{h}$, $\mf{H}$, \sr{hrord}{h}
	    \item $\mf{h}_r^*$, $\mf{H}_r^*$,  $\mf{h}^*$, $\mf{H}^*$, \sr{Hr*}{h*}
	    \item $\mf{h}_{\mf{m}}$, $\h$, \pr{hm}
	    \item $\mf{h}_{\chi}$, \pr{hchi}
	\end{ilist2}
	\item Hecke operators and maps
	\begin{ilist2} 
	    \item $s_r$, $w_{Np^r}$, $\langle j \rangle_r$, \sr{sr}{<j>r}
	    \item $e_r$, \pr{er}
	    \item $U_t$, \pr{Ut}
	    \item $\langle j \rangle$, $\langle j \rangle^*$, \sr{<j>}{<j>*}
	    \item $\Pi$, \pr{Pi}
	\end{ilist2}
	\item homology classes
	\begin{ilist2} 
	    \item $\binom{a}{b}_r$, $\{\alpha,\beta\}_r$, \sr{binom}{path}
	    \item $[u:v]_r$, $[u:v]_r^{\pm}$, \sr{Manin}{Maninpm}
	    \item $\xi_r(u:v)$, $\xi(u:v)$, \sr{xir}{xi}
	    \item $\bar{\xi}(u:v)$, \pr{barxi}
	\end{ilist2}
	\item integers
	\begin{ilist2} 
	    \item $N$, \pr{N}, \pr{N2}
	    \item $p$, \pr{p}, \pr{p2}, \pr{p3}
	\end{ilist2}
	\item Kummer maps
	\begin{ilist2} 
	    \item $\pi_a$, \pr{pia}
	    \item $\Theta_K$, $\phi_2$, \sr{ThetaK}{phi2}
	\end{ilist2}
	\item $L$-functions
	\begin{ilist2} 
	    \item $\mc{L}_N$, $\mc{L}_{N,M}$, \sr{LN}{LNM}
	    \item $\mc{L}_{N,M}^{\star}$, \pr{LNMstar}
	    \item $\overline{\mc{L}_{N,M}^{\star}}$, \pr{overlineLNMstar}
	    \item $L_{p,M}$, $L_{p,M}^{\star}$, $L_p$, \sr{LpM}{Lp}
	    \item $\overline{L_{p,M}^{\star}}$, $\overline{L_{p,M}}$, $\overline{L_p}$, 
			\sr{overlineLpMstar}{overlineLp}
	\end{ilist2}
	\item measures
	\begin{ilist2} 
	    \item $\lambda_N$, \pr{lambdaN}
	    \item $\lambda_{N,M}$, \pr{lambdaNM}, \pr{lambdaNM2}
	\end{ilist2}
	\item modular curves and cusps
	\begin{ilist2} 
	    \item $Y_1^r(N)$, $X_1^r(N)$, $C_1^r(N)$, \sr{Y1r}{C1r}
	\end{ilist2}
	\item modular representation and related modules
	\begin{ilist2}
	 	\item $Z_{\theta}$, $Z'_{\theta}$, $\rho_{\theta}$, $B_{\theta}$, $C_{\theta}$, \sr{Zpt}{Bt}
	\end{ilist2}
	\item ordinary (co)homology
	\begin{ilist2} 	
	    \item $H_1(X_1^r(N);\zp)^{\mr{ord}}$, \pr{ordhomology}	   
	    \item $H_1(N)$, $\mc{H}_1(N)$, \pr{H_1}
	    \item $H^1(X_1^r(N);\zp)^{\mr{ord}}$, \pr{ordcohomology}
             \item $H^1(N)$, $\mc{H}^1(N)$, \pr{H^1}
             \item $H^1_{\et}(N)$, $\mc{H}^1_{\et}(N)$, \pr{H^1et}
	    \item $H_1^{\et}(N)$, $\mc{H}_1^{\et}(N)$, \pr{H_1et} 
	    \item $\mc{Y}_N$, $\mc{Z}_N$, \pr{ZN}
	    \item $Y_r$, \pr{Yr}
	    \item $Y_{\chi}$, \pr{Ychi}
	\end{ilist2}
	\item $p$-adic objects
	\begin{ilist2} 
	    \item $\Z_{p,N}$, $\Lambda_N$, \sr{ZpN}{LambdaN}
	    \item $[j]$, $[j]_r$, \sr{[j]}{[j]r}, \sr{[j]2}{[j]r2}
	    \item $\Lambda_N^{\mf{h}}$, \pr{LambdaNh}
	    \item $\Lambda_N^{\ltheta}$, $\mf{L}_N^{\ltheta}$, \sr{LambdaNt}{LNt}
	    \item $c_N$, \pr{cN}
	    \item $\Z_{p,N}^{\star}$, $\Lambda_N^{\star}$, \sr{LambdaNstar}{ZpNstar}
	\end{ilist2}
	\item Poincar\'e duality pairings
	\begin{ilist2} 
	    \item $(\,\cdot\,,\,\cdot\,)_r$, $\langle\,\cdot\,,\,\cdot\, \rangle_N$, 
			\sr{Poincarer}{PoincareN}
	\end{ilist2}
	\item reciprocity maps
	\begin{ilist2} 
	    \item $\Psi_K$, \pr{PsiK}
	    \item $\Psi_K^{\circ}$, \pr{PsiKcirc}
	\end{ilist2}
	\item rings of character values
	\begin{ilist2} 
	    \item $R_{\chi}$, \pr{Rchi}
	    \item $\mc{O}_N$, \pr{ON}
	\end{ilist2}
	\item specialization maps
	\begin{ilist2} 
	    \item $\tilde{\chi}$, \pr{tildechi}, \pr{tildechi2}, \pr{tildechi3}
	\end{ilist2}
	\item set of primes
	\begin{ilist2} 
	    \item $S$, \pr{S}
	\end{ilist2}
	\item twists
	\begin{ilist2} 
	    \item 
	    $\mc{O}_N(\chi)$, $A(\chi)$, $A'(\chi)$, 
	    \sr{ONchi}{A'chi}	
	\end{ilist2}
	\item unit groups
	\begin{ilist2} 
	    \item $\mc{O}_{E,S}^{\times}$, $\mc{E}_E$, $\mc{U}_K$, \sr{OES}{UK}
	    \item $\mc{C}_K$, \pr{CK}
	\end{ilist2}
	\item units
	\begin{ilist2} 
	    \item $\zeta_d$, $\zeta$, \sr{zetad}{zeta}
	    \item $1-\zeta^v$, \pr{1zeta}
	    \item $\eta_{M,r,t}^{\psi}$,  $\alpha_{r,t}^{Q,\psi}$,  $\eta_{M,t}^{\psi}$,
	    $\alpha_t^{Q,\psi}$, $\alpha_t^{\psi}$, 
			\sr{etaMrt}{alphat}
	\end{ilist2}
\end{ilist}
	
\onecolumn

\renewcommand{\baselinestretch}{1}

\end{document}